\documentclass[12 pt]{article}

\usepackage{amsmath,amsthm,verbatim,amssymb,amsfonts,amscd,diagrams}
\theoremstyle{plain}
\newtheorem{theorem}{Theorem}[section]
\newtheorem{corollary}[theorem]{Corollary}
\newtheorem{lemma}[theorem]{Lemma}
\newtheorem{proposition}[theorem]{Proposition}

\theoremstyle{definition}
\newtheorem{definition}[theorem]{Definition}

\theoremstyle{remark}

\newarrow{ul}---->
\newarrow{Backwards}<----

\newcommand{\td}[1]{\tilde{#1}}
\newcommand{\into}{\hookrightarrow}

\newcommand{\R}{\mathbb{R}}
\newcommand{\bd}{\partial}

\newcommand{\mc}[1]{\mathcal{#1}}

\newcommand{\vg}{\varGamma}

\newcommand{\map}{\operatorname{map}}
\newcommand{\hl}{\operatorname{holink}}
\newcommand{\rhl}{\operatorname{rholink}}
\newcommand{\wt}{\widetilde}

\newcommand{\mf}{\mathfrak}

\begin{document}

\title{Stratified fibrations and the intersection homology of the regular neighborhoods of bottom strata}\label{S: nobundle}
\author{Greg Friedman\\Yale University\\Dept. of Mathematics\\10 Hillhouse Ave\\PO Box 208283\\New Haven, CT 06520\\friedman@math.yale.edu\\Tel. 203-432-6473  Fax:  203-432-7316}
\date{March 4, 2003}
\maketitle

\begin{abstract}
In this paper, we develop Leray-Serre-type spectral sequences to compute the intersection homology of the regular neighborhood and deleted regular neighborhood of the bottom stratum of a stratified PL-pseudomanifold. The $E^2$ terms of the spectral sequences are given by the homology of the bottom stratum with a local coefficient system whose stalks consist of the intersection homology modules of the link of this stratum (or the cone on this link). In the course of this program, we establish the properties of stratified fibrations  over unfiltered base spaces and of their mapping cylinders.  We also prove a folk theorem concerning the stratum-preserving homotopy invariance of intersection homology. 

\end{abstract}



\tableofcontents

\section{Introduction}

The purpose of this paper is to prove the existence of a Leray-Serre-type spectral sequence for the intersection homology (with local coefficients) of the regular neighborhood (and deleted regular neighborhood) of the bottom stratum of a stratified PL-pseudomanifold.  Although it is well known that the regular neighborhood of a space $X$ is homotopy equivalent to $X$, this homotopy equivalence will not necessarily preserve the stratification of a filtered space and so it is not possible to determine the intersection homology of the neighborhood from that of $X$ alone as would be possible for ordinary homology. This reflects, of course, that intersection homology is a finer invariant than homology and designed to suit stratified phenomena. On the other hand, the neighborhood does retain stratified versions of certain other essential properties of regular neighborhoods. For instance, it is homotopy equivalent to the mapping cylinder of a fibration in a manner that preserves the essential properties of the stratification.  The particular  fibration involved can be given in terms of the homotopy links of Quinn \cite{Q1} and provides a special case of the stratified fibrations of Hughes \cite{Hug}, whose properties we must first develop. In fact, we will show that any such regular neighborhood is stratum-preserving homotopy equivalent to the mapping cylinder of such a stratified fibration and then employ the stratum-preserving homotopy invariance of intersection homology. This last seems to be a well-known fact, but its proof appears to be unavailable in the literature and is thus presented below in full generality, including local coefficient systems. 

Since some of the spaces involved in carrying out our program do not appear to be even topological pseudomanifolds (in the sense of \cite{GM2}), it is not clear that there exists a proof of our main theorem (at least along these lines) that utilizes the sheaf-theoretic machinery of intersection homology, so instead we use singular intersection homology as developed in \cite{Ki}; note that on a stratified pseudomanifold, the two theories are isomorphic if we assume sheaf homology with compact supports or work on a compact space. Once we have established the proper stratified analogues of fibration theory and homotopy invariance of intersection homology, many of the subsequent details directly mirror the standard proof of the existence of  spectral sequences for the homology of ordinary fibrations, though there are still several minor modification necessary in order to handle the nuances of intersection homology and local coefficients. 

Our main theorem then is the following; see below for further explanations of the notations:

\begin{theorem}[Theorem \ref{P: Neigh IH}]
Let $X$ be a finite-dimensional stratified pseudomanifold with
locally 
finite triangulation and filtration $\emptyset=X_{-1}\subset X_0\subset\cdots\subset X_n=X$ such that $X_i=\emptyset$ for $i<k$. Let $N=N(X_k)$ be an
open
regular neighborhood of $X_k$, and let $L$ be the link of the stratum $X_k$ (if $X_k$ is not connected, then we can treat each component separately and each component will have its own link). Then, for any fixed perversity $\bar p$ and local coefficient system $\mc{G}$ defined on $X-X_{n-2}$,
there are homological-type spectral sequences $\bar E^r_{p,q}$ and
$E^r_{p,q}$ that abut (up to isomorphism)
to
$IH^{\bar p}_i(N-X_k;\mc{G})$ and $IH^{\bar p}_i(N;\mc{G})$ with respective
$E^2$
terms
\begin{align*}
\bar  E^2_{p,q}=H_p(X_k; \mc{IH}^{\bar p}_q(L;\mc{G}|L)) &&
E^2_{p,q}=H_p(X_k; \mc{IH}^{\bar p}_q(cL;\mc{G}|cL))
\end{align*}
($cL=$ the open cone on $L$),
where $\mc{IH}^{\bar p}_q(L;\mc{G}|L)$ and $\mc{IH}^{\bar p}_q(cL;\mc{G}|cL)$ are local
coefficient systems with respective stalks $IH^{\bar p}_q (L;\mc{G}|L)$ and $IH^{\bar p}_q
(cL;\mc{G}|cL)$ (see below). Furthermore, the map
$i_*: IH^{\bar p}_i(N-X_k;\mc{G})\to
IH^{\bar p}_i(N;\mc{G})$ induced by inclusion induces a map of spectral sequences which
on the the $E^2$ terms is determined by the coefficient homomorphism
$\mc{IH}^{\bar p}_q(L;\mc{G}|L)\to \mc{IH}^{\bar p}_q(cL;\mc{G}|cL)$ given by the map on the
stalk intersection homology modules induced by the inclusion $L\into cL$.
\end{theorem}

The steps to the proof of this theorem proceed as follows: In Section \ref{S: IH inv}, we recall the definition of intersection homology on a filtered space and show that intersection homology is a stratum-preserving-homotopy invariant. In Section \ref{S: strat fib}, we introduce a special case of
the stratified fibrations of Hughes \cite{Hug} and examine some
of its properties. In particular, we show that these satisfy various stratum-preserving analogues of the standard properties of fibrations. We also show, in Section \ref{S: bundles}, how a stratified fibration can determine a bundle of coefficients on its base space with coefficient stalks given by the
intersection homology of the stratified fiber. In Section \ref{S: spaces}, we present some notations
and facts from the works of Quinn \cite{Q1}, \cite{Q2} and Hughes
\cite{Hug} on  manifold weakly stratified spaces and homotopy links. Finally, in Section \ref{S: theorem}, we put this
all together to prove the theorem. 

We have also provided an appendix which contains a corrected version of the proof of the theorem (see  \cite{Ch79} and \cite{Q1}) that certain  deformation retract neighborhoods are homotopy equivalent to the mapping cylinders of certain path spaces.

The spectral sequences developed here can be viewed as fundamental tools for the computation of intersection homology modules, not just for regular neighborhoods but for entire spaces. Coupled with the intersection homology analogue of the Mayer-Vietoris theorem, these spectral sequences may be utilized to carry out computations by working up through a stratified pseudomanifold one stratum at a time. We refer the reader to \cite{GBF2} for an example of such techniques as applied to study the properties of
intersection Alexander polynomials of PL-knots which are not locally-flat. 

Parts of this work originally appeared as part of the author's dissertation. I thank my advisor, Sylvain Cappell, for all of
his generous and invaluable guidance.

For references to  intersection homology theory, the reader is advised to
consult \cite{GM1}, \cite{GM2}, \cite{Ki}, and \cite{Bo}.

\section{Stra\-tum-pre\-serving homotopy
invariance of intersection homology}\label{S: IH inv} 

In this section, we review some of the relevant definitions of intersection homology on a filtered space (note that we do not yet assume the spaces to be stratified in the sense of pseudomanifolds with normally locally trivial strata (e.g. \cite{GM2})). Then we fill in some of the details, which are left
to the reader in \cite{Q2}, on the invariance of intersection homology under
stratum-preserving homotopy equivalence. We also show how this works for intersection
homology with local coefficients.

We first recall some standard definitions. As defined in \cite{GM1}, a (traditional) \emph{perversity} $\bar p$ is a sequence of integers $\{\bar p_0, \bar p_1, \bar p_2, \ldots\}$ such that $\bar p_i\leq \bar p_{i+1}\leq \bar p_i +1$ and such that $\bar p_0=\bar p_1=\bar p_2=0$ (much of the following would hold using more general perversities (see \cite{Ki}, \cite{CS}, \cite{GBF}), but for the other common case, in which $\bar p_2=1$ (superperversities), it is not clear that intersection homology with local coefficients can be well defined in a geometric manner, i.e. without sheaves; hence we omit discussion). 

A \emph{filtered space} is a space, $X$, together with a collection of
closed subspaces
\begin{equation*}
\emptyset = X_{-1}\subset X_0\subset X_1 \subset \cdots \subset X_{n-1}\subset X_n=X.
\end{equation*}
If we want to emphasize both the space and the filtration, we will refer to the
filtered space $(X,\{X_i\})$. Note that $X_i=X_{i+1}$ is possible. We will
refer to $n$ as the (stratified) \emph{dimension} of $X$ and $X_{n-k}$ as the \emph{$n-k$
skeleton} or the
\emph{codimension $k$ skeleton}. The sets $X_i-X_{i-1}$ are the \emph{strata} of $X$. We call a space either \emph{unfiltered} or \emph{unstratified} if we do not wish to consider any filtration on it. 

We can
now define the singular intersection homology of $X$ for a perversity $\bar
p$, $IH_*^{\bar p}(X)$, in the usual manner: as the homology of the chain complex $IC_*^{\bar
p}(X)$, which is the submodule generated by the $\bar p$-allowable chains of the singular chain complex, $C_*(X)$. A singular $i$-simplex
$\sigma:\Delta^i\to X$ is \emph{$\bar p$-allowable} if
$\sigma^{-1}(X_{n-k}-X_{n-k-1})$ is contained in the $i-k+\bar p(k)$ skeleton
of the (polyhedral) simplex $\Delta^i$ (with the usual filtration by polyhedral
skeletons), and an $i$-chain is \emph{$\bar p$-allowable} if it and its boundary (in
the chain complex sense) are linear combinations of $\bar p$-allowable
$i$-simplices and $i-1$ simplices, respectively. Similarly, we can define the
intersection homology with coefficients in a group or with local coefficients (see below for more details). If $X$ is a
stratified PL-pseudomanifold with filtration determined by the
stratification, then this
definition of intersection homology is the standard one (with compact supports, in the sheaf language) and is a topological invariant \cite{Ki}.
Note, however, that for a general filtered space, $X$, the intersection
homology may not be a topological invariant, e.g. it may depend on our choice
of filtration. 

Since the intersection homology definitions depend only on the
codimensions of the strata, we see that we can reindex a filtration by
addition of a  fixed integer, accompanied by the same change to the
stratified dimension, without affecting the intersection homology of the
space. In what follows, we shall often defer from the norm by \emph{not}
reindexing in certain situations where it is standard to do so when working
with stratified pseudomanifolds. In particular, for simplicity we will
usually give
filtered
subspaces, $Y\subset X$, the same dimension as $X$ and filtration indexing
$Y_i=Y\cap X_i$, even though the common practice for stratified
pseudomanifolds is
to reindex to the simplicial dimensions of $Y$ and $Y_i$ by subtracting the
codimension of $Y$ in $X$ from each index. 

If $X$ and $Y$ are two filtered spaces, we call a map $f:X\to Y$
\emph{stratum-preserving} if the image of each component of a stratum of $X$
lies in a stratum of $Y$ (compare \cite{Q1}). In general, it is not required
that strata of $X$ map to strata of $Y$ of the same (co)dimension, but unless
otherwise noted we will assume that $X$ and $Y$ have the same stratified
dimension and that $f(X_i-X_{i-1})\subset Y_i-Y_{i-1}$. We call $f$
a \emph{stratum-preserving homotopy equivalence} if there is a stratum-preserving
map $g:Y\to X$ such that $fg$ and $gf$ are stratum-preserving homotopic to the
identity (where the filtration of $X\times I$ is given by the collection
$\{X_i\times I\}$). We will sometimes denote the stratum-preserving homotopy
equivalence of $X$ and $Y$ by $X\sim_{sphe}Y$ and say that $X$ and $Y$ are
\emph{stratum-preserving homotopy equivalent} or \emph{s.p.h.e.}

A slightly more general notion which we will need below is the following \cite{Q2}:
Consider two filtered spaces $X$ and $Y$ with a common subset $U$. A map
$f:X\to Y$ of filtered spaces is a \emph{stratum-preserving homotopy equivalence near U}
provided:
\begin{enumerate}
\item $f$ is stratum-preserving,
\item $f|U$ is the identity,
\item There exists a neighborhood $V\subset Y$ of $U$ and a stratum-preserving map
$g:V\to X$ such that $fg$ and $gf$, where defined, are stratum-preserving homotopic
rel $U$ to the inclusions (i.e. $fg$ is stratum-preserving homotopic rel $U$ to the inclusion $V\into Y$
and $gf$ is stratum-preserving homotopic rel $U$ to the inclusion of $f^{-1}(V)$ into
$X$).
\end{enumerate}

We now study in some detail the effect of a stratum-preserving homotopy
equivalence on intersection homology. If $X \sim_{sphe} Y$, 
then it follows that there must be a bijective correspondence between
non-empty strata of $X$ and non-empty strata of $Y$ determined by the mapping. 
We prove the following proposition which is implied in \cite{Q2}:

\begin{proposition}\label{P: s.p.h.e. of IH}
If $f: X\to Y$ is a stratum-preserving homotopy equivalence with inverse
$g:Y\to X$, $X$ and $Y$ have the same dimension $n$ (as filtered spaces), and
$f(X_i-X_{i-1})\subset Y_{i}-Y_{i-1}$ and $g(Y_i-Y_{i-1})\subset X_{i}-X_{i-1}$
for each $i$, then $f$ and $g$ induce isomorphisms of intersection
homology
$f_*: IH_i^{\bar p}(X)\to IH_i^{\bar p}(Y)$ and $g_*: IH_i^{\bar p}(Y)\to
IH_i^{\bar p}(X)$.
\end{proposition}
\begin{proof}
Note that, since $X$ and $Y$ are the disjoint unions of their
strata, the assumptions that $f(X_i-X_{i-1})\subset Y_{i}-Y_{i-1}$ and
$g(Y_i-Y_{i-1})\subset
X_{i}-X_{i-1}$
for each $i$ imply that $f^{-1}(Y_i-Y_{i-1})=X_i-X_{i-1}$ for each non-empty
stratum and
similarly for $g$.

We observe that $f_*$ and $g_*$ are well-defined maps on intersection
homology: If $\sigma:\Delta^i\to X$ is a $\bar p$-allowable singular
$i$-simplex of $X$, then $f_*\sigma=f\sigma:\Delta^i \to Y$ is a singular
$i$-simplex of $Y$. To see that it is allowable, we note that if
$Y_{n-k}-Y_{n-k-1}$ is a non-empty stratum of $Y$, then
$(f\sigma)^{-1}(Y_{n-k}-Y_{n-k-1})=\sigma^{-1}f^{-1}(Y_{n-k}-Y_{n-k-1})$. But
$f^{-1}(Y_{n-k}-Y_{n-k-1})= X_{n-k}-X_{n-k-1}$, so
$(f\sigma)^{-1}(Y_{n-k}-Y_{n-k-1})=\sigma^{-1}(X_{n-k}-X_{n-k-1})$, which lies
in the $i-k+\bar p(k)$ skeleton of $\Delta^i$ since $\sigma$ is $\bar
p$-allowable in $X$. Hence $f$ induces a map $f_*: IC_i^{\bar p}(X) \to
IC_i^{\bar p} (Y)$ which clearly also commutes with boundary maps so that it
takes cycles to cycles and boundaries to boundaries, inducing a well-defined
map on intersection homology.

We next show that $(fg)_*$ and $(gf)_*$ are isomorphisms, which will suffice to
complete the proof. By assumption, there exists a stratum-preserving homotopy
$H: X\times I \to X$ such that, for $x\in X$, $H(x,0)$ is the identity and
$H(x,1)=gf(x)$. Now consider a $\bar p$-allowable $i$-cycle $\sigma=\sum_j
a_j\sigma_j$ in $X$, where the $a_j$ are coefficients and the $\sigma_j$ are
$i$-simplices. For the (polyhedral) $i$-simplex $\Delta^i$, let $[v_{0},\ldots,
v_{i}]$ and $[w_{0},\ldots, w_{i}]$ be the simplices $\Delta^i\times 0$ and
$\Delta^i\times 1$ in $\Delta^i\times I$. Then $\Delta^i\times I$ can be
triangulated by the $i+1$ simplices $\Delta^{i+1,l}=[v_{0},\cdots, v_{l},
w_{l}, \cdots, w_{i}]$. Now consider the singular $i+1$ chain in $X$ given
by
\begin{equation} 
F=\sum_j a_j \sum_l (-1)^l H\circ (\sigma_j\times
\text{id})\circ\sigma_{jl}, 
\end{equation} where $(\sigma_j\times
\text{id}):\Delta^i\times I\to X\times I$ is given by $\sigma_j$ on the first
coordinate and the identity on the second and $\sigma_{jl}$ is the singular
$i+1$ simplex given by the inclusion $\Delta^{i+1,l}\into \Delta^i\times I$.
Then, by a simple computation (see \cite[\S 2.1]{Ha} for the details), $\bd
F=(gf)\sigma-\sigma$. Thus $(gf)_*\sigma=\sigma$ in ordinary homology. To show
that this holds in intersection homology, it remains only to check that $F$ is
a $\bar p$-allowable $i+1$ chain. Since $\sigma$ is $\bar p$-allowable by
assumption and $gf\sigma$ is allowable as the image of $\sigma$ under a
stratum-preserving map, $\bd F$ is $\bar p$-allowable, and it suffices to check
that each $H\circ (\sigma_j\times \text{id})\circ\sigma_{jl}: \Delta^{i+1,l}\to
X$ is a $\bar p$-allowable singular simplex, i.e. that
$\sigma_{jl}^{-1}(\sigma_j\times \text{id})^{-1}H^{-1}(X_{n-k}-X_{n-k+1})$ lies
in the
$i+1-k+\bar p(k)$ skeleton of $\Delta^{i+1,l}$.

We must consider the inverse image of each stratum $X_{n-k}-X_{n-k-1}$. Since
we filter $X\times I$ by $\{X_k\times I\}$ and $H$ is a stratum-preserving
homotopy to the identity, $H^{-1}(X_{n-k}-X_{n-k-1})=(X\times I)_{n-k}-(X\times
I)_{n-k-1}=(X_{n-k}-X_{n-k-1})\times I$.  Then $(\sigma_j\times\text{id})^{-1}
[(X_{n-k}-X_{n-k-1})\times I]=[\sigma_j^{-1}(X_{n-k}-X_{n-k-1})]\times I$.
Since $\sigma_j$ is $\bar p$-allowable, this is contained in the product of $I$
with the $i-k+\bar p(k)=:r$ skeleton of $\Delta^i$ (denoted $\Delta^i_r$). Finally, we show that
$\sigma_{jl}^{-1}(\Delta^i_{r}\times I)$ lies in the $r+1$ skeleton of the
$i+1$-simplex $\Delta^{i+1,l}$ by construction. In fact,
$\sigma_{jl}^{-1}(\Delta^i_r\times I)=(\Delta^i_{r}\times I)\cap [v_{0},\ldots,
v_{l}, w_{l}, \ldots, w_{i}]$. But $\Delta^i_r\times I=\cup
[v_{\alpha_0},\ldots, v_{\alpha_r}]\times I$, where the union is taken over all
sets of ordered $r$-tuples of integers with $0\leq \alpha_0<\ldots
<\alpha_r\leq i$. Furthermore, each $[v_{\alpha_0},\ldots, v_{\alpha_r}]\times
I=\cup_{\beta}[v_{\alpha_0},\ldots, v_{\alpha_{\beta}}, w_{\alpha_{\beta}},
\ldots, w_{\alpha_r}]$. That each $[v_{\alpha_0},\ldots, v_{\alpha_r}]\times I$
is a union of such simplices follows as for our original decomposition of
$\Delta^i\times I$, and each such simplex lies in the given triangulation
of $\Delta^i\times I$ since each is clearly a face of one of the
$\Delta^{i+1,l}$.
Therefore, any simplex in the intersection of $\Delta^{i+1,l}$ and
$\Delta^i_r\times I$ is spanned by a set of vertices which is a subset of both
$\{v_{0},\ldots, v_{l}, w_{l}, \ldots, w_{i}\}$ and one of the collections
$\{v_{\alpha_0},\ldots, v_{\alpha_{\beta}}, w_{\alpha_{\beta}}, \ldots,
w_{\alpha_r}\}$. Thus, each such simplex must have dimension less than or equal
to $r+1$, and $(\Delta^i_{r}\times I)\cap \Delta^{i+1,l}$ must lie in the $r+1$
skeleton of $\Delta^{i+1,l}$. Hence $H\circ (\sigma_j\times \text{id})\circ
\sigma_{jl}$ is a $\bar p$-allowable $i+1$-simplex.

This shows that $(gf)_*=g_*f_*$ is an isomorphism of intersection homology (in
fact
the identity). Similarly, $(fg)_*=\text{id}_*$, and so $f_*$ and $g_*$ must each
be isomorphisms.

\end{proof}

This proof clearly extends to hold for any constant coefficient module. To show that this also
works for local systems of coefficients, it is perhaps simplest to work with
an analogue of the following version of ordinary homology with local coefficients as
presented in \cite[\S 3.H]{Ha} (note that Hatcher confines his definition to bundles of abelian groups, but, at least for homology, there is no difficulty extending the definition to bundles of modules):   Given a space $X$ with local coefficient
system $\mc{G}$ with stalk module $G$ and given action of the fundamental group, we
use this data to form a bundle of modules (or ``sheaf space'' in the language
of sheaves), which we also denote $\mc G$. Then an element of the singular
chain complex $C_i(X;\mc{G})$ is a finite linear combination $\sum_j n_j
\sigma_j$, where each $\sigma_j$ is a singular simplex $\sigma_j:\Delta^i\to X$ and
$n_j$ is a lift of $\sigma_j$ to $\mc{G}$. If $n_j$ and $m_j$ are both lifts
of the same singular simplex, then $n_j+m_j$ is defined using the addition on
the bundle, which is continuous. We also take $n_j+m_j=0$ if it is the
lift to the zero section. Similarly, we can define scalar multiplication by elements $r$ in the ground ring $R$ of the bundle of (left) $R$-modules  by defining $n_j\to rn_j$ pointwise and noting the continuity of the scalar multiplication operation in a bundle (clearly similar definitions can be made for bundles of right $R$-modules). The coefficient of a boundary face of a simplex is
given by $\bd (\sum_j n_j \sigma_j)= \sum_{j,k} (-1)^k
n_j\sigma_j|[v_0,\cdots,\hat v_k,\cdots,v_i]$, where each $n_j$ on the right
hand side is the restricted lift $n_j|[v_0,\cdots,\hat v_k,\cdots,v_i]$. With
these definitions, the
homology of the chain complex $C_*(X;\mc{G})$ is the usual homology of $X$
with local coefficients $\mc{G}$. Note that given a map between spaces $f: X\to Y$
covered by a bundle map of coefficient systems $\td f: \mc{G} \to \mc{H}$,
then $f_*(\sum_j n_j \sigma_j)$ is defined as $\sum_j (\td f\circ n_j)
\cdot(f\circ \sigma_j)$ and induces a map on homology $f_*:H_*(X;\mc{G})\to
H_*(Y;\mc{H})$.

Before getting to intersection homology, let us first show the following: 
\begin{proposition} \label{P: loc coeff he}
Suppose we are given homotopy equivalent (unfiltered) spaces $X$ and $Y$,
homotopy inverses $f:X\to Y$ and $g:Y\to X$, and a local coefficient system
$\mc{G}$ over $X$. Then $g$ covered by the induced map $\td
g:g^*\mc{G}\to\mc{G}$ induces a homology isomorphism $g_*:H_*(Y;g^*\mc{G})\to
H_*(X;\mc{G})$.
\end{proposition}
\begin{proof}
Consider the following diagram:
\begin{equation}\label{E: coeff pullbacks}
\begin{CD}
g^*f^*g^*\mc{G} &@>\bar g >> &f^*g^*\mc{G} &@>\td f >> &g^*\mc{G} &@>\td g
>> & \mc{G}\\
@VVV&&@VVV&&@VVV&&@VVV\\  
Y &@>g >> &X &@>f >> &Y &@>g >> & X,
\end{CD}
\end{equation}
where $\td g$, $\td f$, and $\bar g$ are the maps induced by the pullback
constructions. We will first show that the map $gf$ covered by
the map $\td g \td f$ 
induces an isomorphism $(gf)_*:H_*(X;f^*g^*\mc{G})\to H_*(X;\mc{G})$. To do so,
let $H$ denote a homotopy $H: X\times I\to X$ from the identity to $gf$.
Then we have an induced bundle $H^*\mc{G}$ over $X\times I$ and a map $\td
H:H^*\mc{G}\to \mc{G}$ such that, over $X\times 0$, this is simply the
identity map $\mc{G}\to\mc{G}$ and, over $X\times 1$, it is the induced map
$f^*g^*\mc{G}\to \mc{G}$.  

To show the surjectivity of the homology map $(gf)_*=g_*f_*$
induced by $gf$ covered by $\td g\td f$, let $[C]$ be an element of $
H_i(X;\mc{G})$. Then, by our definitions, $[C]$ is represented by a cycle
$C=\sum n_j\sigma_j$, where the $\sigma_j$ are singular simplices, i.e. maps
of $\Delta^i$ into $X$, and the $n_j$ are lifts of the $\sigma_j$ to
$\mc{G}$. We can also then consider these simplices and their lifts as maps
into $X\times 0$ and $H^*\mc{G}|X\times 0$. Next, we can extend each
simplex map $\sigma_j$ to the map $\sigma_j\times \text{id}:\Delta^i\times
I \to X\times I$. Since bundles of coefficients are covering spaces, there
exist unique lifts, $\wt{\sigma_j \times \text{id}}$, of each of these maps
extending the lift over $X\times 0$. Together these provide a homotopy from
$C$ over $X\times 0$ to a new chain, say $\bar C=\sum \bar n_j\sigma_j$,
where these $\sigma_j$ are the same as those above identifying $X$ and
$X\times 1$, and the $\bar n_j$ are given by the lifts $\wt{\sigma_j \times
\text{id}}|\Delta^i\times 1$. Furthermore, $\bar C$ is a cycle: since $C$
is a cycle, the sum of the lifts over each point $x\times 0\in X\times 0$
of $\sum \bd
\sigma_i$ lies in the zero section of the bundle, but, by
the properties of bundles of coefficients, the liftings of the paths
$x\times I\in
X\times I$ collectively provide a homomorphism from $\mc{G}|x\times 0\to
\mc{G}|x\times
1$. So the sum of the lifts over $x\times 1$ at the other end of the 
homotopy is also $0$. Now, we can compose each of the $\sigma_j\times
\text{id}$ covered
by $\wt{\sigma_j \times \text{id}}$ with $H$ covered by $\td H$ to obtain a
homotopy into $(X,\mc{G})$ from the maps representing $C$ to the maps
representing the image
of $\bar C$ under $gf$ covered by $\td g\td f$. By the prism construction  
employed above in the proof of Proposition \ref{P: s.p.h.e. of IH}
to break a homotopy into simplices (and the obvious extension to coefficient
lifts),
this gives a
homology between $[C]$ and the homology class of
the image $g_*f_*(\bar{C})$. Hence $(gf)_*=g_*f_*$ is
surjective.

For injectivity, let $\bar C=\sum \bar n_i\sigma_i$ be a cycle representing an
element of $H_i(X;f^*g^*\mc{G})$ which maps to $0$ in $H_i(X;\mc{G})$ under
$(gf)_*$. As in the last paragraph, we can lift $\sigma_j\times
\text{id}:\Delta^i\times I\to X\times I$ to an extension of $\bar n_i$ in
$H^*\mc{G}$ (although this time we extend in the other direction). Again, this
induces a chain $C$ in $(X\times 0,\mc{G})$, and composing the homotopies and
their lifts with $H$ covered by $\td H$ induces a homotopy and, by the prism
construction and its lift, a homology from the image of $\bar C$ to $C$ (recall
that $H$ and $\td H$ are the identity on and over $X\times 0$, respectively).
But since the image of $\bar C$ is $0$ by assumption, $C$ must bound another
chain with local coefficients, say $D=\sum m_k\tau_k$, where the $\tau_k$ are
singular $i+1$ simplices in $X$ and the $m_k$ are their lifts to $\mc{G}$. Now,
we proceed as in the last paragraph: we can consider $D$ as a chain with
coefficients in $X\times 0$ covered by $\mc{G}$ and lift the maps
$\tau_k\times\text{id}:\Delta^{i+1}\times I\to X \times I$ to maps into
$H^*\mc{G}$
extending the lifts $m_k$. Again this induces a chain $\bar D=\sum \bar m_i
\tau_i$ in $(X\times 1, f^*g^*\mc{G})$, and, due to the unique path lifting
property of covering spaces, it is readily verified that $\bd \bar D=\bar C$.
Thus, $\bar C$ represents the $0$ element of $H_i(X;f^*g^*\mc{G})$, and it follows that  
the homology map $(gf)_*$ is injective.

We have thus shown that the map  $(gf)_*=g_*f_*$ on homology induced by $gf$
covered by
$\td g\td f$ is an isomorphism. By the same reasoning, the map $fg$ covered
by $\td f\bar g$ (see diagram \eqref{E: coeff pullbacks}) induces a
homology isomorphism $(f\bar g)_*=f_*\bar g_*: H_*(Y;g^*f^*g^*\mc{G})\to
H_*(Y;g^*\mc{G})$. But
together these imply that the homology map induced by $f$ covered by $\td
f$ is injective and surjective and hence an isomorphism
$H_*(X; f^*g^*\mc{G})\to H_*(Y;g^*\mc{G})$. Finally, since $f$ covered by  
$\td f$ induces an isomorphism and $gf$ covered by $\td g\td f$ induces an
isomorphism, $g$ covered by $\td g$ must induce an isomorphism.
\end{proof}

For intersection homology with local coefficients (and traditional
perversities), we can apply the same ideas. As in the stratified
pseudomanifold
case, we assume that the singular set has codimension at least two (i.e.
$X_{n-1}=X_{n-2}$) and that a local coefficient system $\mc{G}$ is given on
$X-X_{n-2}$. We can also think of $\mc G$ as a covering space over $X-X_{n-2}$.
Now, given a $\bar p$-allowable $i$ simplex, $\sigma$, it will no longer be
possible, in general, to lift the entire simplex, as $\sigma$ may intersect
$X_{n-2}$. However, we can choose for coefficients lifts of $\sigma\cap
(X-X_{n-2})$, which consists of
$\sigma$ minus, at most, pieces of its $i-2$ skeleton. So, in particular, we
can lift at least the interior of the simplex and the interior of any
 facet ($i-1$ face).
The coefficient of the boundary $\bd \sigma$ can again be defined in terms of
the signed (i.e. plus or minus)
restrictions of the lift of the simplex to its boundary pieces, again with the
limitation that
some lower-dimensional pieces (and only lower-dimensional pieces) might be
missing. In this manner, the definition of the local coefficient of a simplex
by its lift is well-defined for $\bar p$-allowable intersection simplices, and
it is further well-defined (modulo
lower  
dimensional skeleta) on all boundary pieces by taking restrictions
multiplied by $\pm 1$ since 
the allowability conditions ensure that the interior of each facet of an
allowable simplex lies in the top stratum.
Defining
the chain complex $IC_*^{\bar p}(X;\mc{G})$ as the $\bar p$-allowable
chains formed by linear combinations of
singular
simplices with such lifts as coefficients yields a chain complex whose
homology is the intersection homology of $X$ with coefficients in 
$\mc{G}$. This approach is clearly
equivalent to the more common definition of intersection homology with local-coefficient which involves the lifts of
given points of the simplices; it has the disadvantage of the lift not being
defined for all points (though it is for enough!) and the advantage of not
having to nitpick about which points and how to change points under boundary
maps.

We should, however, note the following concerning the definitions of the last
paragraph: When considering $\bd \sigma$ of a $\bar p$-allowable $i$-chain
$\sigma=\sum n_j\sigma_j$ with local coefficients as defined above, a given
boundary piece of a simplex $\sigma_j$ may not itself be allowable, just as is
the case with constant coefficients. However, by restriction and multiplication
by $\pm 1$, it will have a coefficient lifting defined over its intersection
with the top stratum of $X$, which includes at least its interior. But, again
as in the constant coefficient case, all of the non-$\bar p$-allowable pieces
must cancel in $\bd \sigma$ since each non-zero simplex of $\bd \sigma$ must be
$\bar p$-allowable since $\sigma$ was $\bar p$-allowable. In other words, the
coefficients, where defined, over each point of each of the
non-allowable simplices in $\bd \sigma$ must add up to zero. Since the
remaining boundary pieces are
$\bar p$-allowable, their coefficient lifts are again defined, at the least,
over their
interiors and the interiors of their  facets. Therefore, $\bd\bd\sigma$
is well-defined. It is then routine to check that $\bd\bd\sigma=0$ just as for
ordinary homology with or without coefficients, $\bd\bd\sigma$ being equal to
the boundary of what remains of $\bd \sigma$ after the pieces which cancel have
been removed (i.e. set to zero). Alternatively, we could note that $\bd\bd
n_i\sigma_i$ is zero where the lift is defined and apply linearity (where it
is not well-defined, it must be zero in $\bd \sigma$ anyway, by its
non-allowability), or we could
appeal directly to the equivalence of this definition of intersection chains
with coefficients and the more standard definition involving lifts over
well-chosen points.

Suppose  now that the filtered space $X$ is given a local coefficient system
$\mc{G}$, as above, and that $f:X \to Y$ is a stratum-preserving homotopy
equivalence with inverse $g$. We also again assume that $X$ and $Y$ have
the same stratified dimension and that $f$ and $g$ take strata to strata of the
same (co)dimension.
Since $f$ and $g$ induce homotopy equivalences
of $X-X_{n-2}$ and $Y-Y_{n-2}$ (again assuming $X_{n-1}=X_{n-2}$ and
$Y_{n-1}=Y_{n-2}$), they also induce isomorphisms of the fundamental groups and
determine local coefficients on $Y$ via the pullback $g^*$. Furthermore, it is
clear that $f^*g^*\mc{G}$ is bundle isomorphic to $G$.

We want to show that $g$ covered by $\td g:g^*\mc{G}\to\mc{G}$ over the top
stratum induces an isomorphism of intersection homology with local
coefficients. But this now follows from a combination of the techniques from
the proofs of Propositions \ref{P: s.p.h.e. of IH} and \ref{P: loc coeff he}.
In particular, we simply apply the proof of Proposition \ref{P: loc coeff he},
but
all bundles and lifts are restricted to lie over the top stratum. Once again,
we will not necessarily be able to lift entire simplices, but the homology
theory works out as described in the last few paragraphs. Furthermore,
since all of the homotopies in that proof can be taken to be
stratum-preserving in this context, there is no
trouble with extending lifts: any point that is mapped to a top stratum remains
in the top stratum under all of the necessary homotopies and so its lift can be
uniquely extended. Also, there is no difficulty showing that the
homotopies induce homologies by breaking the homotopies
of simplices up into triangulated prisms. The allowability issues are taken
care of just as in the proof of Proposition \ref{P: s.p.h.e. of IH}, and
therefore each prism piece is allowable and its coefficient is determined over,
at the least, its interior and the interiors of its facets as the restriction of the lift
over the entire homotopy on the top stratum.
Therefore, we have shown the following:

\begin{corollary}\label{C: s.p.h.e. with l.c.}
Suppose that $g: Y\to X$ is a stratum-preserving homotopy equivalence, $\mc{G}$
is a local coefficient system on $X-X_{n-2}$ (which is the top stratum), $X$ and $Y$ have the same
dimension $n$ (as filtered spaces), and 
$f(X_i-X_{i-1})\subset Y_{i}-Y_{i-1}$ and $g(Y_i-Y_{i-1})\subset X_{i}-X_{i-1}$
for each $i$. Then the intersection homology map
$g_*: IH_i^{\bar p}(Y ;g^*\mc{G})\to IH_i^{\bar p}(X ;\mc{G})$, induced by $g$
covered over $X-X_{n-2}$ by the natural map $g^*\mc{G}\to \mc{G}$, is an
isomorphism.
\end{corollary}

Although we have have focused on stratum-preserving homotopy equivalences,
which will be our primary use of the above theory, these same arguments can be
easily generalized to prove the following corollary.

\begin{corollary}\label{C: l.c.h} 
If $f_0, f_1: X\to Y$ are stratum-preserving homotopic by a homotopy $H$,
$\mc{G}$ is a local coefficient system over $X-X_{n-2}$, $\mc{G}'$ is a local
coefficient system over $Y-Y_{n-2}$, $X$ and $Y$ have the same stratified
dimension, $f_0(X_i-X_{i-1})\subset Y_{i}-Y_{i-1}$ (or, equivalently,
$f_1(X_i-X_{i-1})\subset Y_{i}-Y_{i-1}$) for each $i$,
and $f_0, f_1$ are covered over $X-X_{n-2}$ by maps $\td
f_0, \td f_1:\mc{G}\to \mc{G}'$ which are homotopic over $H|_{(X-X_{n-2})\times
I}$ by a homotopy $\td
H$, then $f_{0*},f_{1*}: IH_i^{\bar p}(X ;\mc{G})\to IH_i^{\bar p}(Y ;\mc{G}')$
are identical maps on intersection homology.
\end{corollary}
\begin{proof} [Sketch of proof] 
For simplicity, let us omit the direct mention of the fact that the covers only
lie over the top stratum; it will be assumed. Once again, we can employ a prism
so that if $\sigma:\Delta^i\to X$ is a singular cycle with coefficient lift
$\td \sigma$,
then $H\circ (\sigma\times \text{id})$ covered by $\td H\circ (\td \sigma\times
\text{id})$ provides a homotopy between $f_0 \sigma$ with lift $\td f_0\td
\sigma$ and $f_1\sigma$ with lift $\td f_1\td \sigma$. Hence, we can
again employ a prism as in the proof of Proposition \ref{P: s.p.h.e. of IH}
to show that $f_0$ and $f_1$ induce identical maps on homology with local
coefficients. 
\end{proof}

\section{Stratified fibrations}\label{S: strat fib} In this section, we study a special case of the notion of a
stratified
fibration as introduced by Hughes \cite{Hug} (our definition will be less general than that of Hughes, which assumes that both of the spaces in a stratified fibration are filtered, but it will suit our purposes).

\begin{definition}(Hughes, \cite[p. 355]{Hug})
If $Y$ is a filtered space, $Z$ and $A$ unfiltered spaces, a map $f: Z\times A \to Y$ is \emph{stratum-preserving along $A$} if, for each $z\in Z$, $f(z\times A)$ lies in a single stratum of $Y$.
\end{definition}

\begin{definition}
Suppose $Y$ is a filtered space and $B$ is an unfiltered space (equivalently
a space
with one stratum). A map $p:Y\to B$ is a \emph{stratified fibration} if, given any space
$Z$ and the commutative diagram
\begin{equation}\label{E: strat fib}
\begin{CD}
Z &@>f>> &Y\\
@V\times 0 VV && @VVpV\\
Z\times I & @>F>> &B,
\end{CD}
\end{equation}
there exists a \emph{stratified solution} which is stratum-preserving along $I$, i.e. a map $\td F: Z\times I\to Y$ such that $p\td F=F$ and, for each $z\in Z$, $\td F(z,0)=f(z)$ and  $\td F( z\times I)$ lies in a single stratum of $Y$.
\end{definition}

We want to show that stratified fibrations satisfy stratified analogues of several of the standard properties of fibrations, such as the existence of pullbacks and triviality over a contractible base. 
We obtain these results largely by reproducing the proofs for standard fibrations found in Whitehead \cite{Wh}, adding the relevant details where the stratifications come into play.

First, we present the following two lemmas which are special cases of \cite[Lemmas 5.2 and 5.3]{Hug} of Hughes and follow immediately from them:

\begin{lemma}
Given two solutions, $\td F$ and $\td G$, of a stratified lifting problem
\eqref{E: strat fib}, then there exists a homotopy $H:Z\times I\times I\to Y$
from $\td F$ to $\td G$ rel $Z\times\{0\}$ such that $pH(z,t,s)=F(z,t)$
and $H$ is stratum-preserving along $I\times I$.
\end{lemma}

\begin{lemma}[Stratified relative lifting]\label{L: extend}
Suppose $p:Y\to B$ is a stratified fibration, $B$ unstratified, $(Z,A)$ an unstratified NDR-pair, and we are given the commutative diagram:
\begin{equation*}
\begin{CD}
Z\times \{0\}\cup A\times I &@>f >> & Y\\
@VVV && @VVpV\\
Z\times I &@>F >> & B.\\
\end{CD}
\end{equation*}
If there is a deformation retraction $H_t: Z\times I\times I\to Z\times I$ of 
$Z\times I$ to
$Z\times \{0\}\cup A\times I$ such that $fH_1:Z\times I \to Y$ is
stratum-preserving along $I$, then there exists a stratified solution $\td
F:Z\times I\to Y$ which is stratum-preserving along $I$, satisfies $p\td F=F$,
and $\td F| Z\times \{0\}\cup A\times I=f$. (Note: Hughes assumes in his lemma
that $Z$ is a metric space, but I don't believe this is ever used in the
proof, at least supposing the definition of NDR as given in \cite{Wh}.)
\end{lemma}

In fact, Hughes' proof of Lemma \ref{L: extend} applies to the following slightly more general situation: 

\begin{lemma}\label{L: DR ext} Suppose $p:Y\to B$ is a stratified fibration, $(Z\times I,A)$ is a DR-pair, $i:A\to Z\times I$ is the inclusion, and we are given the commutative diagram:
\begin{equation*}
\begin{CD}
A &@>f >> & Y\\
@ViVV && @VVpV\\
Z\times I &@>F >> & B.\\
\end{CD}
\end{equation*}
If there is a deformation retraction $H_t: Z\times I\times I\to Z$ of $Z\times
I$ to $A$ such
that $fH_1:Z \times I\to Y$ is stratum-preserving along $I$, then there exists a
stratified solution $\td F:Z\times I\to Y$ which is stratum-preserving along $I$, satisfies $p\td F=F$, and $\td F| A=f$.
\end{lemma}

We now define the induced stratified fibration in analogy with induced
fibrations. 

\begin{definition} Suppose $p:Y\to B$ is a stratified fibration and $f:B'\to B$
is a map
of spaces. We denote by $Y'$ (or $f^*Y$ if we want to emphasize the map $f$)
the set $\{(b',y)\in B'\times Y|f(b')=p(y)\}$ and call it the \emph{induced 
stratified fibration} or \emph{pullback stratified fibration}. We filter $Y'$
as a subset of
$B'\times Y$, i.e. $Y'_i=\{(b',y)\in Y'| y\in Y_i\}$. 
\end{definition}
We obtain maps $f':Y'\to
Y$ and $p':Y'\to B'$ induced by the projections. Note that $f'$ is obviously
stratum-preserving and takes the fiber over $b'\in B'$ to the fiber over $f(b')$. 
If $f$ is an inclusion, then so is $f'$.
Assuming this notation, we prove the following basic facts:

\begin{lemma}[Universality] If $g:Z\to Y$ is a stratum-preserving map and $q:Z\to B'$ is
a map such that $pg=fq$, then there is a unique stratum preserving map $h:Z\to Y'$ such that $f'h=g$ and $p'h=q$. 
\end{lemma}
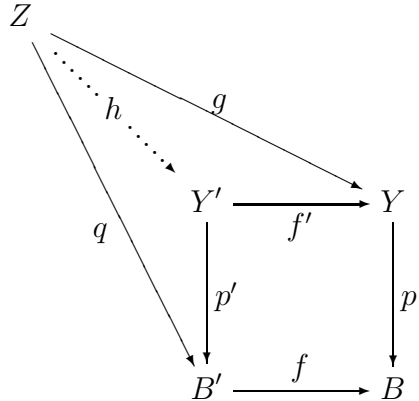
\begin{proof}

\begin{diagram}
Z \\
  &\relax\rdDotsto~{h}     \rdTo(4,2)^g
    \rdTo(2,4)_q  \\
  &         & Y' & \rTo_{f'} & Y     \\
  &         & \dTo_{p'}             &        & \dTo_p \\
  &         & B'                 & \rTo^f & B      \\
\end{diagram}
Based on the hypotheses, $q$ and $g$ define a map $(q,g):Z\to B'\times Y$ which
is stratum-preserving with respect to the product filtration on $B\times Y$.
Clearly, Im$(q,y)\subset Y'$, and thus $(q,g)$ induces a stratum-preserving map
to $Y'$. The commutativity relations are obvious. Uniqueness follows obviously
from $Y'$ being a subset of $Y\times B'$: the two projections $f'$ and $p'$
define the image point uniquely in $Y'$, but these are fixed by the hypotheses.
\end{proof}

\begin{lemma}
The induced map $p':Y'\to B'$ is a stratified fibration.
\end{lemma}
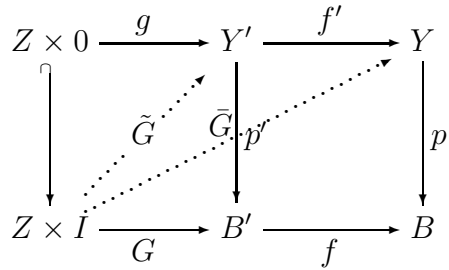
\begin{proof}
\begin{diagram}
Z\times 0 &\rTo^g & Y' &\rTo^{f'} & Y\\
\dInto &\ruDotsto(2,2)~{\td G}&\dTo_{p'} & \ruDotsto(4,2)^{\bar G} & \dTo_p\\ 
Z\times I &\rTo_G & B' &\rTo_f & B\\
\end{diagram}

If $g:Z\to Y'$ and $G:Z\times I\to B'$ are the data for a stratified
lifting problem for the map $p'$, then the maps $fG$ and $f'g$ are the
data for a stratified lifting problem for $p$, which has a stratified
solution $\bar G$ by assumption. Now, stratify $Z\times I$ by $(Z\times
I)_i=\bar G^{-1}(Y_i)$. Note that if $(z,t)\in (Z\times I)_i$ for some $z\in
Z$
and $t\in I$, then $z\times I\in (Z\times I)_i$ by the definition of a
stratified lifting. With this stratification, $\bar G$ is clearly a
stratum-preserving map. But then $\bar G:Z\times I\to Y$ and $G:Z\times
I\to
B'$ satisfy $fG=p\bar G$, the hypothesis of the Universality Lemma. Hence
there is a unique stratum-preserving map $\td G:Z\times I\to Y'$ such that
$p'\td G=G$ and $f'\td G=\bar G$, and, since $\bar G|Z\times \{0\}= f'g$ and
$G|Z\times 0=p'g$, $\td
G|Z\times\{0\}\to Y'$ is given by $\td G(z\times
\{0\})=(p'\td G(z\times 0),f'\td G(z\times 0))=(p'g(z\times\{0\}), f'g(z\times
\{0\}))\in Y'\subset B'\times Y$.
But since a point in $B'\times Y$ is determined uniquely by its
projections, this is also $g(z\times\{0\})$, so $\td G$ extends $g$. Finally,
since $\td G$ is
stratum-preserving on the strata we have imposed on $Z\times I$, each $z\times I$,
lying in a single stratum of $Z\times I$, must map into a single stratum of
$Y'$
under $\td G$. Hence, $\td G$ is a stratified solution to the lifting
problem and $p'$ is a stratified fibration.
\end{proof}

\begin{corollary} The restriction of a stratified fibration is a stratified fibration. \end{corollary}
\begin{proof}
Pull back over the inclusion which induces the restriction.
\end{proof}
\begin{corollary}\label{cor: trans}
Suppose that $p:Y\to B$ is a stratified fibration, $f:B'\to B$ and $g:B''\to
B'$ are maps, $p':Y'\to B'$ is the stratified fibration induced by $p$ and $f$,
and $p'':Y''\to B''$ is the stratified fibration induced by $p'$ and $g$. Then
up to natural stratified isomorphism, $p''$ is also the stratified fibration
induced from $p$ by the map $fg$.
\end{corollary}
\begin{proof}
As for ordinary fibrations, $\phi: (fg)^*Y\to g^*f^*Y$ given by
$\phi(b_2,y)=(b_2,(g(b_2),y))$ provides the isomorphism. It is clearly
stratum-preserving.
\end{proof}

\begin{definition}
Suppose that $p_i: Y_i\to B$, $i=1,2$, are stratified fibrations (note that these subscripts do not indicate strata). A map
$H: Y_1\times I\to Y_2$ is a \emph{stratum-preserving fiber homotopy}
(\emph{s.p.-homotopy}) if it is
stratum-preserving and $p_2H(y,t)=p_1(y)$ for all $t\in I$. $f:Y_1\to Y_2$ is a
\emph{stratum-preserving fiber homotopy equivalence} if $f$ is a
stratum-preserving map such that $p_2f(y_1)=p_1(y_1)$ for all $y_1\in Y_1$ and there is a
stratum-preserving map $g:Y_2\to Y_1$ such
that $p_1g(y_2)=p_2(y_2)$ for all $y_2\in Y_2$ and $fg$ and $gf$ are stratum-preserving fiber homotopic
to the respective identity maps on $Y_2$ and $Y_1$. 
\end{definition}

\begin{lemma}
Suppose that $f_0, f_1:B'\to B$ are homotopic maps, $p:Y\to B$ is a stratified
fibration, and let $p'_t:Y'_t\to B'$, $t=0,1$ be the induced pullback
stratified
fibrations. Then $Y'_0$ and $Y'_1$ are stratum-preserving fiber homotopy
equivalent, i.e. there exist stratum-preserving fiber homotopy
equivalences $g: Y'_0\to Y'_1$ and $h:Y'_1\to Y'_0$ and stratum-preserving
fiber homotopies between $gh$ and $hg$ and the respective identity
maps. In this situation, we will say that $p'_0$ and $p'_1$ have the same
\emph{stratified fiber homotopy type} or are \emph{stratum-preserving fiber homotopy
equivalent (s.p.f.h.e.)}.
\end{lemma}
\begin{proof}
The proof follows \cite[I.7.25]{Wh} closely.
Let $i_t:B'\to B'\times I$, $t=0,1$, be the inclusions of the ends, and let
$f:B'\times I\to B$ be a homotopy from $f_0$ to $f_1$. Then $fi_t=f_t$.  Using Corollary \ref{cor: trans}, it then suffices to prove that for any stratified fibration $p:Y\to B\times I$, if $p_t:Y_t\to B$ are the stratified fibrations induced by $i_t:B\to B\times I$ and $Y_t$ has the filtration induced as a subset of $Y$, then $p_0$ and $p_1$ have the same stratified fiber homotopy type. So suppose that this is the set-up.

Let $i'_t$ be the inclusion maps $Y_t\to Y$ over the inclusions $i_t$. The maps
$\text{id}\times p_0:I\times Y_0\to I\times B$ and the inclusions $f_0:Y_0\to
Y$ form
the data for a stratified lifting problem, and, by assumption, there exists a
stratified solution $H_0: I\times Y_0\to Y$ such that $H_0(0,y_0)=i'_0(y_0)$
and $pH_0(s,y_0)=(s, p_0(y_0))$. Since $i'_0$ is the inclusion, $H_0$ takes
each $I\times y_0$ into the stratum of $Y$ in which $y_0$ lies. Similarly,
there is an $H_1$ such that $H_1(1,y_1)=i'_1(y_1)$, $pH_1(s,y_1)=(s,
p_1(y_1))$, and $I\times y_1$ maps into the stratum containing $y_1$. Define
$\phi:Y_0\to Y_1$ and $\psi:Y_1\to Y_0$ by $i'_1\phi(y_0)=H_0(1,y_0)$ and
$i'_0\psi(y_1)=H_1(0,y_1)$. From the properties of $H_t$, $\phi$ and $\psi$ are
stratum-preserving, $p_1\phi=p_0$, and $p_0\psi=p_1$.

To see that $\psi\phi$ is stratum-preserving fiber homotopic to the identity,
define $g: 1\times I\times Y_0\cup I\times \dot I\times Y_0\to
Y$ and $K:I\times I\times Y_0\to I\times B$ by $K(s,t,y_0)=(s, p_0(y_0))$ and 
\begin{equation*}
g(s,t,x)=
\begin{cases}
H_0(s, y_0), & t=0\\
H_1(s,\phi(y_0)), & t=1\\
i'_1\phi(y_0), & s=1.
\end{cases}
\end{equation*}
$g$ is well defined and stratum preserving on the induced strata. These maps
form the data for a stratified lifting extension problem.  By Lemma \ref{L:
extend}, it has a solution if there is a deformation retraction $R_t: I\times
I\times Y_0 \times I\to I \times I\times Y_0$ to $1\times I\times Y_0\cup
I\times \dot I\times Y_0 $ such that $gR_1:I\times I\times Y_0\times I \to Y$
is stratum-preserving along $I$. But we can take as $R$ the map $r\times
\text{id}_{Y_0}$ where $r$ is any standard deformation retraction of $I\times
I$ into $1\times I\cup I\times \dot I$. Then each $I\times I\times y_0\times I$
clearly gets mapped into a single stratum of $I\times I\times Y_0$ under $R$
and hence under $R_1$, and $g$ is stratum-preserving from its definition. Hence
the conditions are
met for there to be a stratified extension of $g$, namely $G:I\times I\times
Y_0\to Y$, such that $pG=K$ and $G$ is stratum-preserving along the first $I$
factor. But it is also stratum-preserving along the second $I$ factor since
clearly $G$ is stratum-preserving along it for $s=1$. Then $G|0\times I\times
Y_0$ is a stratum-preserving fiber homotopy between $\psi\phi$
and the identity. The other case, $\phi\psi$, is handled similarly.
\end{proof}

\begin{corollary}
If $p:Y\to B$ is a stratified fibration and $B$ is contractible then $p$ is stratified fiber homotopically trivial, i.e. there is a stratum-preserving fiber homotopy equivalence between $p:Y\to B$ and $p':F\times B\to B$ for some stratified fiber $F$. 
\end{corollary}
\begin{proof}
Clearly the pullback of a map from $B$ to a point in $B$ gives the product stratified fibration.
\end{proof}

We will call a stratum-preserving fiber homotopy equivalence given by the map $\phi:F\times B\to Y$ a
\emph{stratified trivialization}. We call it a \emph{strong stratified trivialization} if, for some
$b_0\in B$ with
$F=p^{-1}(b_0)$, we have $\phi(x,b_0)=x$ for all $x\in F$.

\begin{corollary}\label{C: strong triv}
If $(B, b_0)$ is a DR-pair and $p:Y\to B$ is a stratified fibration, then there exists a strong stratified trivialization $\phi: p^{-1}(b_0)\times B\to Y$. Any two such strong stratified trivializations are stratum-preserving fiber-homotopic rel $p^{-1}(b_0)\times b_0$. 
\end{corollary}
\begin{proof}
The proofs are analogous to those given for ordinary fibrations in \cite[I.7.28-29]{Wh}; the necessary modifications to the stratified case are either obvious or employ the same techniques used in the preceding proofs. 
\end{proof}

Finally, we will need some results on the mapping cylinders of stratified
fibrations. These are the subject of the next few lemmas. 

We will always use mapping cylinders (including cones) with the \emph{teardrop topology} (see \cite[\S 2]{Hug99} or \cite[\S 3]{HTWW}). 
The \emph{teardrop topology} on the mapping cylinder $I_{p}$ of the map $p: Y\to B$ is  the topology generated by the sub-basis consisting of open subsets of $Y\times (0,1]$ and sets of the form $[p^{-1}(U)\times (0,\epsilon)]\cup U$, where $U$ is an open set in $B$. This is, in fact, a basis since if $U$ and $V$ are two open sets of $B$, then within the intersection $(p^{-1}(U)\times (0, \epsilon_1)\cup U)\cap (p^{-1}(V)\times (0, \epsilon_2)\cup V)$, we have the basis element $p^{-1}(U\cap V)\times (0, \text{min}(\epsilon_1, \epsilon_2))\cup (U\cap V)$. The intersection of any basis elements of the first and second types or of two elements of the first type is a basis element of the first type.

If $p:Y\to B$ is a
stratified fibration, $Y$ filtered by the sets $\{Y_i\}$, then the mapping cylinder
$I_p$ is naturally filtered by $B$,
as the bottom stratum, and the mapping cylinders of $p$ restricted to each
$Y_i$. There remains, however, the question of how to label the dimensions. Even if
we create a consistent plan for how to handle the dimension of each $I_{p|Y_i}$, it
is not clear what the dimension of $B$ should be. Fortunately, in the applications
below, it will actually be the mapping cylinder which comes equipped which a natural
filtration of this form, and the filtration of $Y$ can be considered to be induced
by the intersection of that filtration with $Y\times t\subset I_p$, $t\in
(0,1]$. So, although this approach appears somewhat backwards, we will assume in
what follows that the filtrations arise in this manner.

We will also have need to consider a local coefficient system on the top stratum
$I_p-I_{p,n-2}=(Y-Y_{n-2})\times (0,1]$. By standard bundle theory, any such system,
$\mc{G}$, will be isomorphic to $(\mc{G}|Y\times t)\times (0,1]$, for any choice of  $t\in
(0,1]$. Therefore, we will always assume that $\mc{G}$ has the form of such a
product. Conversely, if we start with a coefficient system, $\mc{G}$, on
$Y-Y_{n-2}$, then 
$\mc{G}\times (0,1]$ gives a coefficient system on $I_p$, unique up to
isomorphism. We will denote this coefficient system by $I_{\mc G}$.
If $I_p$ is a closed cone, $\bar cY$, we will sometimes denote this induced coefficient
system by $\bar c\mc{G}$ \label{Gcone} (note that we reserve the symbol $cX$ for the open cone on the space $X$).

\begin{lemma}\label{L: prod/cyl}
Suppose $f:X\to Y$ is a continuous map with mapping cylinder $I_f$ and that $Z$ is another space. Then $I_f\times Z=I_{f\times \text{id}_Z}$, the mapping cylinder of $f\times
\text{id}_Z:X\times Z\to Y\times Z$.
\end{lemma}

\begin{proof}
Both spaces clearly have the same underlying sets, so we focus on the topologies. $I_f\times Z$ has a basis of the form $V\times W$, where $V$ is a basis element of $I_f$ and $W$ is a basis element of $Z$. In particular, these all have the form either $V_1\times W$ or $V_2\times W$, where $V_1$ is open in $X\times (0,1]$ and $V_2$ has the form $[f^{-1}(U)\times (0,\epsilon)]\cup U$ for an open subset $U\subset Y$. On the other hand, $I_{f\times \text{id}_Z}$ has basis elements of the form $A=\{\text{an open subset of }(X\times Z)\times (0,1]\}$ or $[B=(f\times \text{id}_Z)^{-1} (C)]\times (0,\epsilon)\cup C$ for some open $C\subset Y\times Z$. So we need to show that $A$ and $B$ are open in $I_f\times Z$  and $V_1\times W$ and $V_2\times W$ are open in $I_{f\times \text{id}_Z}$.

Since the subset topologies on $X\times (0,1]\times Z\subset I_f\times Z$ and $X\times Z\times (0,1]\subset I_{f\times \text{id}_Z}$ all include the standard basis elements for products, these topologies are finer than (in fact equal to) the obvious product topologies, and so the sets of type $A$ and $V_1\times W$ are open in the desired spaces.

The set $V_2\times W=[f^{-1}(U)\times (0,\epsilon)\cup U]\times W$ can be rewritten as $[(f\times \text{id}_Z)^{-1}(U\times W)\times (0,\epsilon)]\cup (U\times W)$. This is an element of type $B$, so the inclusion of $V_2\times W$ into $I_{f\times \text{id}_Z}$ is open. 

Finally, the sets of the form $U\times W$ for $U$ open in $Y$ and $W$ open in $Z$ form a basis for $Y\times Z$, so we can restrict to looking at sets of type $B$ in which $C$ has the form $U\times W$. But then  $B=[(f\times \text{id}_Z)^{-1} (U\times W)\times (0,\epsilon)]\cup U\times W=[(f^{-1}(U)\times W\times (0,\epsilon)]\cup U\times W=(f^{-1}(U)\times (0,\epsilon)\cup U)\times W$, which is of the form $V_2\times W$. Hence this inclusion is also open.
\end{proof}

\begin{lemma}
Suppose $f:Y_1\to Y_2$ is a stratum- and fiber-preserving map of stratified
fibrations, $p_i: Y_i\to B$, $i=1,2$. Let $\bar f$ be the map from $I_{p_1}$
to
$I_{p_2}$ induced on the mapping cylinders by the identity on $B$ and $f\times
\text{id}_{(0,1]}: Y_1\times (0,1]\to Y_2\times (0,1]$. Then $\bar f$ is a continuous
stratum-preserving map which, furthermore, takes each cone
$I_{p_1|p_1^{-1}(b)}$, for $b\in B$, into the cone $I_{p_2|p_2^{-1}(b)}$.
\end{lemma}
\begin{proof}
The only inobvious part of the statement of the lemma is that $\bar f$ should be continuous at points in $B$. If $b\in B$, then the basis elements of the topology of $I_{p_2}$ containing $f(b)=b$ have the form $V=[p_2^{-1}(U)\times (0,\epsilon)]\cup U$, where $U$ is a neighborhood of $b$ in $B$. Then from the definitions, $\bar f^{-1}(V)=[f^{-1}p_2^{-1}(U)\times (0,\epsilon)]\cup \bar f^{-1}(U)=[p_1^{-1}(U)\times (0,\epsilon)]\cup U$. 
Here we have used that $\bar f|B$ is the identity on $B$ and that $p_2f=p_1$.  So $\bar f^{-1}(V)$ is an open neighborhood of $b$ that maps into $V$, and $\bar f$ is continuous.

\end{proof}

\begin{corollary}
The map $\bar f$, as defined in the preceding lemma, induces a map on intersection homology
$\bar f_*:IH_*^{\bar p}(I_{p_1})\to IH_*^{\bar p}(I_{p_2})$. If $f$ is covered by a
bundle
map on the top stratum $\td f:\mc{G}\to \mc{H}$, then $\td f\times
\text{id}_{(0,1]}$ gives a map $I_{\mc G}\to I_{\mc{H}}$ covering $\bar f$ on the
top stratum, and together these induce a map $f_*: IH_*^{\bar p}(I_{p_1};I_{\mc
G})\to IH_*^{\bar p}(I_{p_2}; I_{\mc{H}})$. 
\end{corollary}

\begin{lemma}\label{C: cyl fib hom eq}
Suppose that $f: Y_1\to Y_2$ and $g:Y_2\to Y_1$ are stratum-preserving fiber
homotopy inverses for the stratified fibrations $p_i:Y_i\to B$, $i=1,2$. Then
$\bar f$ and $\bar g$ as defined in the previous lemma are stratum- and
cone-preserving
homotopy inverses between the $I_{p_i}$, meaning that the homotopies from
$\bar f\bar g$ and $\bar g\bar f$ to the identities are stratum-preserving and
take
$I_{p_i|p_i^{-1}(b)}\times I$ to $I_{p_i|p_i^{-1}(b)}$.
\end{lemma}
\begin{proof}
Let $H:Y_1\times I\to Y_1$ be the stratum-preserving fiber homotopy from
$gf$ to the identity. By Lemma \ref{L: prod/cyl}, $I_{p_1}\times I$ is equal to
the mapping cylinder of $p_1 \times \text{id}_I: Y_1\times I\to 
B\times I$. Define $\bar H: I_{p_1}\times I\to
I_{p_1}$ so that for each $t\in I$, $\bar H|I_{p_1}\times t$ is the map
determined
as in the previous lemma by extending the map $H| Y_1\times t$ to the cylinder. 

Again, continuity is clear except for points $(b,t)\in  B\times I\subset  I_{p_1}\times I$. Consider $\bar H(b,t)=b$. Again, the basis neighborhoods of $b$ in $B$ have the form $V=[p_1^{-1}(U)\times (0,\epsilon)]\cup U$. Since by Lemma \ref{L: prod/cyl}, we can also consider $I_{p_1}\times I$ as the mapping cylinder of $p_1\times \text{id}_I$, $(b,t)$ has a neighborhood of the form $W=[p_1^{-1}(U\times I)\times (0,\epsilon)]\cup (U\times I)$. Since $H$ is a fiber homotopy, $\bar H$ takes $U\times I$ to $U$ by projection, and since $\bar H$ is cone-preserving and takes each point in the cylinder to a point of the same ``height'', it maps points in $p_1^{-1}(U\times I)\times (0,\epsilon)$ to $p_1^{-1}(U)\times (0,\epsilon)$. Hence the neighborhood $W$ maps into the neighborhood $V$, and $\bar H$ is continuous. 

Furthermore, over $ B\times 0$, since the restriction of $H$ is the identity map of $Y_1$, the restriction of $\bar H$ is the identity map of the
mapping
cylinder $I_{p_1}$, while the restriction over $B\times 1$ of $H$ is $gf$ so
that the restriction of $\bar H$ is $\overline{gf}$, which is clearly $\bar 
g\bar f$ from the
definition of these maps. Therefore, $\bar H$ is a homotopy from
$\overline{gf}$ to the
identity. That it is cone-preserving follows immediately from the definitions
and the fact that $H$ was fiber preserving. Similarly, the fact that $H$ was
stratum-preserving and the definition of $\bar H$ shows that $\bar H$ is
stratum-preserving. The proof that $\bar f\bar g$ is stratum- and
cone-preserving
homotopic to the identity on $I_{p_2}$ is the same.

\end{proof}

\begin{corollary}\label{C: cyl strong triv}
If a stratified fibration $p:Y\to B$ possesses a (strong) stratified trivialization $f:F\times B\to Y$, then there
is also a (strong) stratum- and cone-preserving homotopy equivalence $
\bar cF\times B\to I_p$.
\end{corollary}   
\begin{proof}
By the preceding lemma, there exists a stratum- and cone-preserving homotopy
equivalence $\bar f$ from the mapping cylinder of the projection
$F\times
B\to B$ to $I_p$. Using Lemma \ref{L: prod/cyl}, this mapping cylinder is
exactly
$\bar cF\times B$. If the trivialization is strong so that $f$ restricted to the
fiber over some $b_0$ is the identity map $F\times b_0\to F_{b_0}$, then it is
clear from the construction in the lemma that $\bar f|\bar cF_{b_0}$ will be the
identity
map $\bar cF\times b_0\to \bar cF_{b_0}$.
Furthermore, if the homotopies determining the equivalence of $Y$ and $F\times
B$ are stationary over $b_0$, then so will be the homotopies which determine
the equivalence of $I_p$ and $\bar cF\times B$ as constructed in the proof of the
lemma.
\end{proof}

\begin{corollary}\label{C: s.p.f.hom}
If $H:Y_1\times I\to Y_2$ is a stratum-preserving fiber homotopy from $f$ to $g$,
$f,g: Y_1\to Y_2$ stratum- and fiber-preserving maps, then $\bar
H:I_{p_1}\times I\to Y_2$,
defined so that $\bar H|(I_{p_1}\times t)=\overline{H|I_{p_1}\times t}$ for each
$t\in I$, is a stratum-preserving homotopy $\bar H: I_{p_1}\times I\to I_{p_2}$
which maps $\bar c p_1^{-1}(b)\times t$ to $\bar c p_2^{-1}(b)$ for each $b\in B$
and $t\in I$. 
\end{corollary}
\begin{proof}
The method of proof is the same as that employed in the proof of the lemma to
construct the homotopy from $\bar g\bar f$ to the identity. 
\end{proof}

\begin{corollary}\label{C: s.p.f.hom IH}
If $f$ and $g$ are stratum-preserving fiber homotopic maps of stratified fibrations, then $\bar f$ and $\bar g$
induce the same map on intersection homology $f_*=g_*:IH_*^{\bar p}(I_{p_1})\to
IH_*^{\bar p}(I_{p_2})$. If the homotopy, $H$, from $f$ to $g$ is covered on the top stratum by a
homotopy of local coefficients $\td H:\mc{G}\times I\to \mc{H}$ from $\td f$ to
$\td g$, then $\bar H$ can be covered on the top stratum by a homotopy $\td
H\times \text{id}_{(0,1]}: I_{\mc{G}}\times I\to I_{\mc{H}}$ from $\td f\times
\text{id}_{(0,1]}$ to $\td g\times \text{id}_{(0,1]}$.
Then
$\bar f$ covered by $\td f\times
\text{id}_{(0,1]}$ and $\bar g$ covered by $\td g\times \text{id}_{(0,1]}$ induce
the same
map on intersection homology $f_*=g_*:IH_*^{\bar p}(I_{p_1};I_{\mc
G})\to IH_*^{\bar p}(I_{p_2}; I_{\mc{H}})$. 
\end{corollary}
\begin{proof}
This follows from Corollaries \ref{C: l.c.h} and \ref{C: s.p.f.hom}. 
\end{proof}

\section{Bundles of coefficients induced by a stratified fibration}\label{S: bundles} To construct
bundles of modules from a stratified fibration, we will again need to generalize some of
the standard arguments from the case of ordinary fibrations. We again modify some of the
relevant results from \cite{Wh}.  Given a stratified fibration $p:Y\to B$, with $B$
unfiltered, and a path $u:I\to B$ from $b_0$ to $b_1$, let $F_t=p^{-1}(b_t)$ and $i_t:
F_t\into Y$ (obviously these are stratum-preserving maps if we define the
stratifications on the
$F_t$ by their intersections with the strata of $Y$). We call a stratum preserving map
$h: F_1\to F_0$ $u$-s.p. admissible if there exists a stratum-preserving homotopy
$H:I\times F_1\to Y$ from $h$ to $i_1$ such that $pH(t,y)=u(t)$ for all $(t,y)\in
I\times F_1$. Just as for ordinary fibrations, the following are clear:
\begin{enumerate} 
\item If $u$ is the constant path and $h$ is the identity map $F_0\to F_0$, then $h$ is  $u$-s.p. admissible.
\item If $v:[0,1]\to B$ and $u:[1,2]\to B$ are paths $I\to B$ with $u(1)=v(1)$, $h:
F_1\to F_0$ is $v$-s.p. admissible, and $k: F_2\to F_1$ is $u$-s.p admissible, then
$h\circ k:F_2\to F_0$ is $v*u$-s.p. admissible. 
The proof, as usual, consists of
adjoining two homotopies:
If $H$ and $K$ are the homotopies which provide the admissibility of $h$ and $k$,
then the admissibility of $hk$ is given by the homotopy $L:[0,2]\times F_2\to Y$
defined by
\begin{equation*}
L(t,x)=
\begin{cases}
H(t, k(x)), & t\in [0,1]\\
K(t,x), & t\in[1,2].
\end{cases}
\end{equation*}
Since $H$, $K$, and $k$ are stratum-preserving, so is $L$.
\end{enumerate}

Given $u\in \pi_1(B; b_0, b_1)$, there always exists a $u$-sp admissible map: Let $f:I \times p^{-1}(b_1)\to B$ be given by $f(t,y)=u(t)$, and let
$g:1\times p^{-1}(b_1)$ be given by $g(1,y)=y$. Then this defines a stratified
homotopy lifting problem, which has a solution $G:I\times p^{-1}(b_1)\to Y$ by
assumption. Then define $g': p^{-1}(b_1)\to p^{-1}(b_0)$ by $g'(y)=G(0,y)$. This is a
stratum-preserving map.

We next show that if $u_0$ and $u_1$ are homotopic rel $\dot{I}$ and $h_0$ and $h_1$
are respectively $u_0$- and $u_1$-s.p. admissible, then $h_0$ and $h_1$ are
stratum-preserving homotopic: Let $F:I\times I\to B$ be a homotopy from $u_0$ to
$u_1$ rel $\dot I$, and, for $s=0,1$, let $H_s$ be the homotopy from $h_s$ to $i_1$
lying over $u_s$. Define a map $G:I\times I\times F_1\to B$ by $G(s,t,y)=F(s,t)$
and a map $H: (\dot{I}\times I\cup I\times 1)\times F_1\to Y$ by
\begin{equation*}
H(s,t,y)=\
\begin{cases}
H_s(t,y),& s=0,1\\
y, & t=1.
\end{cases}
\end{equation*} 
These provide the data for a stratified lifting extension problem.
Since $p$ is a stratified fibration, there is a stratified solution, $\bar H$, by Lemma \ref{L:
extend}:  We can deformation retract $I\times I\times F_1$ to $(\dot{I}\times
I\cup I\times 1)\times F_1$ as the product of $F_1$ with any deformation retraction
$I\times
I\to \dot{I}\times   
I\cup I\times 1$, and clearly this retraction composed with $H$ is
stratum-preserving. $H|I\times 0\times F_1$ then gives a stratum-preserving
homotopy from $h_0$ to $h_1$.

Note, incidentally, that together the above properties imply that the fibers of $Y$ are
stratum-preserving homotopy equivalent.

The upshot of this discussion is that given a homotopy class of paths rel endpoints
in $B$, we have defined a unique stratum-preserving homotopy class of maps from the fiber
over the terminal point to the fiber over the initial point. Furthermore, this is
functorial with respect to composition of classes of paths. Hence, we have
established a functor from the fundamental groupoid of $B$ to the category of
filtered spaces and stratum-preserving homotopy classes of maps. This, in turn,
induces a functor $\mf{F}$ to abelian groups or modules by taking, for $b\in B$,
$\mf{F}(b)=IH_i^{\bar p}(p^{-1}(b))$ and as map $\mf{F}(u):IH_*^{\bar
p}(p^{-1}(b_1))\to IH_*^{\bar p}(p^{-1}(b_2))$ the map on intersection homology induced by the stratum-preserving homotopy class of maps determined by $u$ (see Corollary \ref{C: l.c.h}). As
defined in \cite[Chapter VI]{Wh}, such a functor
determines a bundle of groups (modules) over $B$ with fiber isomorphic to
$IH_i^{\bar p}(p^{-1}(b))$ for any $b$ in $B$.

Now suppose that $Y=(Y,\{Y_n\})$, itself, is given a local coefficient system,
$\mc{G}$, over $Y_n-Y_{n-2}$. It will again be convenient to think of $\mc{G}$ as a
space
(in fact a covering space of $Y-Y_{n-2}$ with projection $\pi:\mc{G}\to Y$). In the
following, $\mc{G}$ will always lie
over the top stratum, but for convenience we suppress this from the notation. Let
$b_0$, $b_1$ continue to denote points of $B$ and $F_s=p^{-1}(b_s)$, and consider
the restrictions $\mc{G}|F_s$. As above,  given a path $u$ in
$B$ from $b_0$ to $b_1$, we can construct a stratum-preserving map $g':F_1\to F_0$
as the end of a stratum-preserving homotopy $G:F_1\times I\to Y$ from the
inclusion of $F_1$. The homotopy $G\circ(\pi\times \text{id}):\mc{G}|_{F_1}\times I\to
Y$, together with
the inclusion map $\mc{G}|F_1\into \mc{G}$, gives the data for a lifting extension
problem
which has a solution since $\mc{G}$ is a covering space. Hence, there is a map $\td
G: (\mc{G}|F_1)\times I$ which covers $G$ over the top stratum and which is the
inclusion on $(\mc{G}|F_1)\times 1$. Let $\td
g: \mc{G}|F_1\to\mc{G}|F_0$ be given by $\td g(z)=\td G(z,0)$. Note that, by the
properties of bundles of coefficients, $\td g$ is a bundle map and, in fact,
provides an isomorphism from the stalk over each $x\in F_1$ to the stalk over
$g'(x)$. Together, $g'$ and $\td g$ induce a map $g_*: IH_*^{\bar
p}(F_1;\mc{G}|F_1)\to IH_*^{\bar p}(F_0;\mc{G}|F_0)$ (recall that we
filter each $F$ by its intersections with the filtration of $Y$). Furthermore, $g_*$
is independent of the choice of $G$ and even the choice of $u$ in its relative path
class. We have seen that $g'$ is unique up to stratum-preserving homotopy by
producing a homotopy $H': I\times I\times F_1\to Y$ from $G$ to the homotopy
which gives any other admissible map. Hence, given any such stratum-preserving
homotopy, $H'$, from $G$, we can lift it again to a homotopy $I\times I\times \mc{G}|F_1\to \mc{G}$
which covers $H'$ and is the inclusion on each $I\times 1 \times
\mc{G}|F_1$ by the unique path lifting property of covering spaces. On $I\times
0\times\mc{G}|F_1$, this gives a homotopy from $\td g$ covering the
stratum-preserving homotopy from
$g'$ to the other candidate admissible map. By Corollary \ref{C: l.c.h}, $g'$
covered by $\td g$ and the maps at the other end
of the homotopy induce the same map on intersection homology with local
coefficients. Therefore, in particular, each path class $\xi\in \pi_1(B;b_0,b_1)$
induces a unique map $IH^{\bar P}_*(F_1;\mc{G}|F_1)\to IH^{\bar
P}_*(F_0;\mc{G}|F_0)$.

We show functoriality of these intersection homology maps as induced by path classes
in $B$. It is clear, by these constructions, that the map over the constant path
is the identity. If $v:[0,1]\to B$ and $u:[1,2]:\to B$ are paths such that
$v(1)=u(1)$, we prove that the composition of the induced maps $IH_*^{\bar
p}(F_2;\mc{G}|F_2)\to IH_*^{\bar p}(F_1;\mc{G}|F_1)$ and $IH_*^{\bar
p}(F_1;\mc{G}|F_1)\to IH_*^{\bar p}(F_0;\mc{G}|F_0)$ is the same as the map induced by
the path $v*u$: Let $g'_2:F_2\to F_1$, $g'_1:F_1\to F_0$, and $g'_3:F_2\to F_0$ be
the maps induced as above by $v$, $u$, and $v*u$, respectively, as the ends of
respective stratum-preserving homotopies $G_2:F_2\times [1,2]\to Y$, $G_1:F_1\times [0,1]\to Y$,
and $G_3:F_2\times [0,2]\to Y$. We already know that $g'_3$ and $g'_1g'_2$ are
stratum-preserving homotopic since they are each $v*u$ admissible maps. In fact, we
know that $g'_3$ is given as the end of a homotopy $G_3:F_2\times I\to Y$ over
$v*u$, while
$g'_1g'_2$ can be defined as the end of the stratum-preserving homotopy, say
$G_{1,2}$,
which on $[1,2]$ is given by $G_2$ over $u$ and on $[0,1]$ is the composition of
$G_1$
over $v$ with the map $g'_2\times\text{id}:F_2\times [1,2]\to F_1\times
[1,2]$. Since these each determine $v*u$-admissible maps, we have seen above how to
construct a homotopy $\bar G: I\times
[0,2]\times F_2\to Y$ from $G_3$ to $G_{1,2}$  which is the inclusion on each
$I\times 2\times F_2$ and which gives a  stratum-preserving homotopy
between $g'_3$ and $g'_1g'_2$ when restricted to $I\times 0\times F_2$. But by
covering space theory, we can cover $\bar G$ with a homotopy $I\times [0,2]\times
\mc{G}|F_2\to \mc{G}$ which is the inclusion on $I\times 2\times \mc{G}|F_2$, and
since the
path lifting is unique, the map $I\times 0\times \mc{G}|F_2$ provides a homotopy
between the lifts induced over $g_3$ by the lift over $G_3$ and the lifts induced
over $g'_1g'_2$ by lifting $G_{1,2}$ (which is clearly the same as the composition
of
the lifts over $g'_1$ and $g'_2$ given by the lifts of $G_1$ and $G_2$). Hence
$g'_3$ and $g'_1g'_2$, together with their lifts as constructed above, give
stratum-preserving homotopic maps $F_2\to F_0$ covered by homotopic coefficient
maps, and therefore induce the same map on intersection homology with local
coefficients by Corollary \ref{C: l.c.h}. Thus, the composition of the maps on 
intersection homology induced
by $u$ and $v$ is the same as the map induced by $v*u$. This shows that these maps
define a functor
from
the fundamental groupoid to the category of modules and, therefore, a local system
of coefficients on $B$. This generalizes the above construction to include the case
where $Y$ starts with its own local coefficient system (over $Y-Y_{n-2}$).
Note that
these properties imply that the induced maps $IH_i^{\bar
p}(F_1;\mc{G}|F_1)\to IH_i^{\bar
p}(F_0;\mc{G}|F_0)$ are isomorphisms. We will call these coefficient systems
$\mc{IH}_*^{\bar p}(F;\mc{G})$.

If $I_{p}$ is the mapping cylinder of a stratified fibration $Y\to B$, we can extend
these constructions to obtain bundles of coefficients $\mc{IH}_*^{\bar p}(\bar
cF;\mc {G})$ over $B$ corresponding to the projection $\Pi:I_{p}\to B$ as follows:
For each of the stratum-preserving maps $h: F_1\to F_0$ determined by a path in $B$,
let $\bar ch:\bar cF_1\to \bar cF_0$ be the map determined by the coning. This will
also be stratum-preserving assuming the usual induced stratification on a cone.
Furthermore, any stratum-preserving homotopy of $h$ induces a stratum-preserving
homotopy of $\bar ch$ by Corollary \ref{C: s.p.f.hom} (noting that a cone is
trivially a mapping cylinder of a stratified fibration).
Since the maps $h$ were
uniquely determined up to stratum-preserving homotopy, we see that the maps $\bar
ch$ are also so well-defined. Furthermore, the maps and homotopies which cover these
in $\mc{G}$ similarly extend to maps and homotopies which cover these in $\bar
c\mc{G}$ (see page \pageref{Gcone}) in the obvious manner.
Hence they induce unique maps $IH_*^{\bar p}(\bar
cF_1;\bar{c} \mc G)\to IH_*^{\bar p}(\bar cF_0;\bar{c} \mc G)$ by Corollary \ref{C: s.p.f.hom
IH}. Altogether, this once again forms a
bundle of coefficients over $B$, which we denote $\mc{IH}_*^{\bar
p}(\bar c F;\bar c\mc{G})$. Furthermore, there is clearly a bundle map
$\mc{IH}_*^{\bar p}(F;\mc{G})\to\mc{IH}_*^{\bar p}(\bar c F;\bar c\mc{G})$ induced
by the inclusion $Y\into I_p$.

\section{Weakly stratified spaces}\label{S: spaces}
We next introduce some definitions due to Quinn \cite{Q1},
\cite{Q2}.

A map of pairs of spaces $f: (A,B)\to (X,Y)$ is said to be
\emph{strict} if $f(B)\subset Y$ and $f: (A-B)\subset X-Y$. The space of such strict
maps with the compact-open topology is denoted $\map_s(A,B;X,Y)$. We can then define the
\emph{homotopy link} $\hl(X,Y)=\map_s(I,\{0\};X,Y)$, i.e. the space of maps of the
unit interval into $X$ which take the point $0$ into $Y$ and all other points into the
complement of $Y$. If $X$ is a filtered space, we will also be interested in the
\emph{stratified homotopy link} $\hl_s(X,Y)$, which is the subspace of $\hl(X,Y)$
consisting of maps $I\to X$ which take $0$ into $Y$ and $(0,1]$ into a single stratum
of $X$. $\hl_s(X,Y)$ is itself filtered by the subsets which take $(0,1]$ into a
particular
skeleton of $X$. Note that these are subspaces of $X^I$ and so are metrizable
if $X$ is. There are natural projections $\pi:\hl(X,Y)\to Y$ and
$\pi_s:\hl_s(X,Y)\to Y$ defined by evaluation at $0$. Where there is no room for
ambiguity,
we will also denote $\pi_s$ by $\pi$.

A subspace $Y\subset X$ is \emph{tame} if $Y$ has a neighborhood $N\subset X$ such that
there exists a strict map $(N\times I,Y\times I\cup N\times 0)\to (X,Y)$ which is
the
identity on $Y$ and the inclusion on $N\times 1$, i.e. if there is a deformation
retraction $N\to Y$ such that all points in $N-Y$ remain in $X-Y$ until time 0 (such
a deformation retraction is called \emph{nearly strict}). If $X$
is a filtered space, then it is \emph{weakly stratified} if for each $k>i$,
$(X_{i+1}-X_i)\subset (X_{k+1}-X_{k})\cup(X_{i+1}-X_i)$ is tame and the projection
$\pi: \hl((X_{k+1}-X_{k})\cup(X_{i+1}-X_i),X_{i+1}-X_i)\to X_{i+1}-X_i$ is a fibration.
The filtration of a space $X$ is said to be a \emph{manifold filtration} if each stratum
is a manifold (we will always assume without boundary), and  $X$ is called a 
\emph{manifold weakly stratified space} if it is weakly stratified by a manifold
filtration. As observed in \cite[p. 235]{Q2}, any \emph{locally cone-like} space is
manifold weakly stratified and hence, in particular, stratified
PL-pseudomanifolds are
manifold weakly stratified (see below). 

Now according to Hughes \cite[Cor. 6.2]{Hug}, for a weakly stratified metric
space $X$ with a finite number of strata, $\pi:\hl_s(X,X_k)\to X_k$ is a
stratified fibration. Therefore, if $\mc{G}$ is a bundle of coefficients on
$\hl_s(X,X_k)-\hl_s(X,X_k)_{n-2}$, $n$ is the stratified dimension of $\hl_s(X,X_k)$
(which is the same as the stratified dimension of $X$), and
we let $\mc L$ denote a stratified fiber, then we can apply the results of the preceding
sections to obtain bundles of coefficients $\mc{IH}_*^{\bar p}(\mc L;\mc G)$ and
$\mc{IH}_*^{\bar p}(\bar c\mc L;\bar c\mc G)$ over $X_k$. These are the bundles of
coefficients of the statement of Theorem \ref{P: Neigh IH}, below.

\section{The intersection homology of the regular neighborhood of
the
bottom stratum of a stratified PL-pseudo\-manifold} \label{S: theorem} 

Let us recall some more definitions. Let $c(Z)$ denote the open cone on the space $Z$,
and let $c(\emptyset)$ be a point. 
Generally, a \emph{stratified paracompact Hausdorff space}
$Y$ (see \cite{GM2} or \cite{CS}) is defined
by a filtration
\begin{equation*}
Y=Y_n\supset Y_{n-1} \supset Y_{n-2}\supset \cdots \supset Y_0\supset Y_{-1}=\emptyset
\end{equation*}
such that for each point $y\in Y_i-Y_{i-1}$, there exists a distinguished neighborhood
$N$ of $y$ such there is a compact Hausdorff space $G$ (called the \emph{link} of the component of the stratum), a filtration 
\begin{equation*}
G=G_{n-i-1}\supset  \cdots \supset G_0\supset G_{-1}=\emptyset,
\end{equation*}
and a homeomorphism
\begin{equation*}
\phi: \R^i\times c(G)\to N
\end{equation*}
that takes $\R^i\times c(G_{j-1})$ onto $Y_{i+j}$. 

A \emph{PL-pseudomanifold} of dimension $n$ is a PL space $X$ (equipped with a class of locally finite triangulations) containing a closed PL subspace $\Sigma$ of codimension at least 2 such that $X-\Sigma$ is a PL manifold of dimension $n$ dense in $X$. A \emph{stratified PL-pseudomanifold} of dimension $n$ is a PL pseudomanifold equipped with a filtration such that $\Sigma=X_{n-2}$ and the local normal triviality conditions of a stratified space hold with the trivializing homeomorphisms $\phi$ being PL homeomorphisms. In fact, for any PL-pseudomanifold $X$, such a stratification always exists such that the filtration refines
the standard filtration of $X$ by $k$-skeletons with respect to some triangulation \cite[Chapter I]{Bo}. Furthermore, intersection homology is known to be a topological invariant of such spaces; in particular, it is invariant under choice of triangulation or stratification (see \cite{GM2}, \cite{Bo}, \cite{Ki}).

Now, suppose that the
filtered space
$(X,\{X_i\})$ is a stratified PL-pseudomanifold with each $X_i$ a
subpolyhedron, so that, in
particular, there exists a triangulation of $X$ in which each $X_i$ is
triangulated as a subcomplex. Note that $X$ is also a filtered polyhedron
in the sense considered by Stone in \cite{Sto}. Given such a filtered
polyhedron and
another subpolyhedron $Y$, then the filtered subpolyhedron $(V,\{V_i\})$
is defined (as in Stone \cite[p. 4]{Sto}) to be a \emph{regular
neighborhood}
with \emph{frontier} fr$(V)$ if there are a triangulation $B$ of $X$ in
which $Y$ and the $X_i$ are triangulated by full subcomplexes $A$ and
$B_i$ and a first derived subdivision $B'$ in which
\begin{enumerate}
\item $V=|\bar N(Y,B')|$,
\item $V_i=|\bar N(Y,B'_i)|$,
\item $\text{fr}(V)=|\text{lk}(Y,B')|$,
\end{enumerate} 
where, as usual, $|Z|$ represents the topological space underlying the
simplicial complex $Z$ and $\bar N(Y,B')$ is the (closed) first derived
neighborhood of $Y$. (Stone in fact defines a more general relative
regular neighborhood, but we will not need the added
generality.) It is further shown in \cite{Sto} that such neighborhoods
satisfy the standard properties for regular neighborhoods, such as
uniqueness via PL-isotopies and the generalized annulus property. We can
also define an open regular neighborhood  by taking the
interior of $V$ (meaning here $V-\text{fr}(V)$) and its intersection with
each of the $V_i$. In what follows, ``open regular neighborhood'' will
always mean an open regular neighborhood in this sense.

Finally, let $X$ be a finite-dimensional stratified PL-pseudomanifold with
locally 
finite triangulation and such that $X_i=\emptyset$ for $i<k$ (in
particular $X\subset \mathbb{R}^N$, for some $N$, and hence $X$ is
metrizable). Then $X_k$ is
both the lowest dimensional non-empty skeleton and stratum, and furthermore,
by the definition of a stratified PL-pseudomanifold, it must be a manifold
without
boundary. ($X_k$ must be an open manifold by definition, but it must also be
a
closed subset of $X$.) We will assume that $X_k$ is
connected for notational convenience, but what follows could be applied to each component
separately. As a stratified PL-pseudomanifold, the stratified dimension of
$X$ is
its dimension as a simplicial complex, and similarly for the dimensions of
the filtering subsets, $X_i$.
Let $N=N(X_k)=(N,\{N_i\})$ be an
open
regular neighborhood of $X_k$ (in the sense of the last paragraph). Such a neighborhood exists since a closed regular neighborhood
 $(\bar N, \{\bar N_i\})$ exists by an application of
\cite[Prop. 1.11]{Sto} (in fact, in Stone's system, we would relabel our
$X_k$ as $X_0$, and then the desired neighborhood would be the ``normal
$(n,0)$-regular neighborhood system'' of this proposition).  We also let $L$
denote the link of the stratum $X_k$ so that every point
in $X_k$ has an
open neighborhood PL-homeomorphic to $D^k\times cL$, where $D^k$ is the open
$k$-dimensional ball and $cL$ is the open cone on $L$. We wish to calculate
the intersection homology of $N$ and of $N-X_k$. If $N$ has the topological
structure of a fiber bundle over $X_k$ with fiber $cL$ then this can be computed by sheaf theoretic methods using spectral sequences (these calculations appear in Sections 5.7 and 5.8 of the author's dissertation \cite{GBF}). We now show that such spectral sequences exist without the need for such an assumption.

\begin{theorem}\label{P: Neigh IH}
With $N$, $X$, and $L$ as above and for any fixed perversity $\bar p$, which
we sometimes omit from
the notation, and local coefficient system $\mc{G}$ defined on $X-X_{n-2}$,
there are homological-type spectral sequences $\bar E^r_{p,q}$ and
$E^r_{p,q}$ which abut (up to isomorphism)
to
$IH^{\bar p}_i(N-X_k;\mc{G})$ and $IH^{\bar p}_i(N;\mc{G})$ with respective
$E^2$
terms
\begin{align*}
\bar  E^2_{p,q}=H_p(X_k; \mc{IH}^{\bar p}_q(L;\mc{G}|L)) &&
E^2_{p,q}=H_p(X_k; \mc{IH}^{\bar p}_q(cL;\mc{G}|cL)),
\end{align*}
where $\mc{IH}^{\bar p}_q(L;\mc{G}|L)$ and $\mc{IH}^{\bar p}_q(cL;\mc{G}|cL)$ are local
coefficient systems on $X_k$ with respective stalks $IH^{\bar p}_q (L;\mc{G}|L)$ and $IH^{\bar p}_q
(cL;\mc{G}|cL)$. Furthermore, the map
$i_*: IH^{\bar p}_i(N-X_k;\mc{G})\to
IH^{\bar p}_i(N;\mc{G})$ induced by inclusion induces a map of spectral sequences which
on the the $E^2$ terms is determined by the coefficient homomorphism
$\mc{IH}^{\bar p}_q(L;\mc{G}|L)\to \mc{IH}^{\bar p}_q(cL;\mc{G}|cL)$ given by the map on the
stalk intersection homology modules induced by the inclusion $L\into cL$.
\end{theorem}

\begin{proof} 
Due to the invariance of intersection homology under
stratum-preserving homotopy equivalences, we can replace $N$ by any
stratum-preserving homotopy equivalent space equipped with the coefficients
pulled back by the equivalence map. For simplicity of notation below, we will
also refer to this bundle as $\mc{G}$ where it is unnecessary to make a
distinction. Which bundle is meant should be clear, however, from the space
under consideration. To obtain a space stratum-preserving homotopy equivalent (s.p.h.e.) to $N$, we employ some of the
constructions of Quinn. In particular, by \cite[3.2]{Q1}, since $N$ is
weakly stratified (in fact, locally cone-like), there is a neighborhood $V$ of
$X_k$ in $N$ and a nearly stratum-preserving deformation retraction of $V$ to
$X_k$. This is a deformation retraction given by a homotopy $H:V\times I\to
X$ such that $H|X_k\times I$ is the projection to $X_k$, $H|V\times 1$ is
the inclusion $V\into X$, $H|V\times 0$ maps $V\times 0$ into $X_k$,
and $H|V\times (0,1]$ is stratum preserving along $(0,1]$. In other words,
the retraction
is stratum preserving until the last moment at time $0$ when everything must
collapse into $X_k$. We will use this to prove the following lemma:

\begin{lemma}\label{l: def ret}
There is a nearly stratum-preserving deformation retraction $R$ of $N$ to $X_k$ in 
$N$, i.e. $R: N\times I\to N$.
\end{lemma}

Assuming the lemma for now, we can obtain a stratum-preserving homotopy equivalence between $N$ and the mapping
cylinder, $I_{\pi}$, of the stratified fibration $\pi: \hl_s(N,X_k)\to X_k$. The idea of this equivalence is found in work by Quinn \cite{Q1} and Chapman \cite{Ch79}. We present in the appendix to this paper a proof suggested by Bruce Hughes that fixes certain technical problems in this earlier work, specifically the lack of continuity of the holink evaluation $I_{\pi}\times I\to N$, $(w,s)\to w(s)$. 
Let us here at least describe stratum-preserving homotopy inverses $f:N\to I_{\pi}$ and $g:I_{\pi}\to N$. Given a point $y\in N$, let $R_y$ be the retraction path of $y$ under the nearly stratum-preserving deformation retraction $R$. For $x\in N$, let $d(x,X_k)$ be the distance from $x$ to $X_k$; we are free
to choose our regular neighborhood (or metric) so that this is always less
than $1$.  Define $f(y)$ by $f(y)=y$ if $y\in X_k$ and $f(y)=(R_y, d(y, X_k))$ if $y\in N-X_k$ (recall that we use the teardrop topology on all mapping cylinders). In order to define $g$, we use an intermediate map $S: \hl_s(N,X_k)\times (0,1]\to \hl_s(N,X_k)$. This homotopy should be thought of as shrinking the paths in the holink so that at time $1$ we have the identity and as we approach time $0$ the paths shrink to nothing. The shrinking is done in a certain way as to make the map $g$ continuous; see the appendix for details. Define $g$ by $g(y)=y$ for $y\in X_k$ and for $(w,s)\in \hl_s(N_k,X)\times (0,1]$ by $g(w,s)=S(w,s)(s)$; this is the evaluation of the path $S(w,s)$ at the time $s$.

Let the stratified dimension of $I_{\pi}$ be the dimension of $X$,
and filter $I_{\pi}$ by $X_k$ and the
mapping cylinders of $\pi|\hl_s(N,X_k)_i$ (note each such $i$ is greater
than $k$). With these definitions, $f(N_i-N_{i-1})\subset
I_{\pi,i}-I_{\pi,i-1}$, and similarly for $g$.
Then, by the proof in the appendix, 
$g:I_{\pi} \to N$ and $f: N\to I_{\pi}$ are stratum preserving homotopy
inverses rel $X_k$. In particular then, by Corollary
\ref{C: s.p.h.e. with l.c.},
$g$ must induce an
isomorphism of intersection homology modules 
with the local coefficient system on $I_{\pi}$ given
by $g^*\mc{G}$. As noted above, for simplicity of notation we will also
denote $g^*\mc{G}$ by  
$\mc{G}$ depending on the space under discussion to indicate which is
meant. Note also that $g^*\mc{G}\cong (g^*\mc{G}|\hl_s(N,X_k))\times (0,1]$
by a unique isomorphism, so we will assume that it has this form. Since these
homotopy
equivalences are stratum-preserving, their restrictions induce a stratum
preserving homotopy-equivalence between $N-X_k$ and $I_{\pi}-X_k$, which is
also clearly stratum-preserving homotopy equivalent to $\hl_s(N,X_k)$. Hence,
again by Corollary \ref{C: s.p.h.e. with l.c.}, there are isomorphisms
$IH_*^{\bar p}(\hl_s(N,X_k);\mc{G})\cong IH_*^{\bar
p}(I_{\pi}-X_k;\mc{G})\cong IH_*^{\bar p}(N-X_k;\mc{G})$, the former induced
by inclusion, the latter by $g$.

Letting $\Pi$ stand for the obvious projection $I_{\pi}\to X_k$, we now
filter $I_{\pi}$ by open sets. It
is unfortunate that we deal with two distinct concepts commonly referred to
as filtrations. As our purpose is to construct spectral sequences, we will
briefly refer to this filtration as an \emph{SS-filtration}. We will also
use the term
\emph{simplicial
skeleton} to refer to the skeleton of a polyhedron in the sense of skeleta of
simplicial complexes. Consider the triangulation of $X_k$, and let $B^s$ be
the regular neighborhood of the simplicial $s$-skeleton of $X_k$ given by its
star on the second derived subdivision of the triangulation. These provide
an SS-filtration of $X_k$ and induces an SS-filtration on $I_{\pi}$ by
$J^s=\Pi^{-1}(B^s)$:
\begin{equation*}
\begin{diagram}
I_{\pi}=J & \supset \cdots \supset & J^s &\supset & J^{s-1} & \supset\cdots\supset & J^0
&
\supset & \emptyset \\
\dTo &   & \dTo & & \dTo & & \dTo &
& \dTo \\
X_k=B & \supset \cdots \supset & B^s &\supset & B^{s-1} & \supset\cdots\supset & B^0
&
\supset & \emptyset. \\
\end{diagram}
\end{equation*}
These are both SS-filtrations by open sets, and thus we can define the singular
intersection chain complexes
$IC_i(J^s;\mc{G})$ and filter $IC_i(I_{\pi};\mc{G})$ by
$F_sIC_*(I_{\pi};\mc{G})=\text{im}(IC_*(J^s;\mc{G})\to IC_*(I_{\pi};\mc{G}))$. This is
obviously an
increasing and exhaustive SS-filtration. The general theory of spectral sequences of
SS-filtrations (see \cite{Mc}) then provides a first quadrant spectral sequence with
$E^1$
term
\begin{multline*}
E^1_{p,q}=H_{p+q}(F_pIC_*(I_{\pi};\mc{G})/F_{p-1}IC_*(I_{\pi};\mc{G}))\\=
H_{p+q}(IC_*(J^p,J^{p-1});\mc{G})=IH_{p+q}(J^p,J^{p-1};\mc{G})
\end{multline*}
and boundary map, $d^1$, given by the boundary map of the
long exact intersection homology
sequence of the triple $(J^s,J^{s-1},J^{s-2})$. (Note that the long exact sequence of the
triple for intersection homology can be defined just as in the usual case by using
quotients of the intersection chain complexes.) The spectral sequence abuts to
$H_*(IC_*(I_{\pi};\mc{G}))=IH_*(I_{\pi};\mc{G})$. Similarly, we can filter
$\hl_s(N,X_k)$
by $\bar J^s=\pi^{-1}B^s$ and in the same manner obtain a spectral sequence
abutting to
\begin{equation*}
H_*(IC_*(\hl_s(N,X_k);\mc{G}))=IH_*(\hl_s(N,X_k);\mc{G})\cong IH_*(I_{\pi}-X_k;\mc{G})
\end{equation*}
and with $\bar E^1_{p,q}=IH_{p+q}(\bar J^p,\bar J^{p-1};\mc{G})$.

We can view each $\bar J^s$
as the subset $\bar J^s\times 1$ of the mapping cylinder
$I_{\pi|\pi^{-1}(B^s)}$, and its filtration is
the same as the induced filtration as a subset. The same is true for any pair 
$\pi^{-1}(U)$ and $\Pi^{-1}(U)$, for $U$ an open subset of $B$. Hence these
inclusions induce well-defined maps on intersection homology and, in
particular, a map of spectral sequences which abuts to the map
$IH_*(\hl_s(N,X_k);\mc{G})\to IH_*(I_{\pi}; \mc{G})$ which, employing the
homotopy equivalences and their induced isomorphisms, is isomorphic to the
map $IH_*(N-X_k;\mc{G})\to IH_* (N;\mc{G})$
induced by inclusion. On the $E^1$ terms, this map is that induced by the
inclusions $IH_{p+q}(\bar J^p,\bar J^{p-1};\mc{G})\to
IH_{p+q}(J^p,J^{p-1};\mc{G})$.

It remains to identify the $E^2$ and $\bar E^2$ terms of these sequences and the 
induced map $\bar E^2\to E^2$. In
particular, we need to show that $E^1_{*,q}$ and $\bar E^1_{*,q}$ are isomorphic to
the chain complexes with local coefficients
$\text{Cell}_*(X_k;\mc{IH}_{q}(cL;c\mc{G}))$ and 
$\text{Cell}_*(X_k;\mc{IH}_{q}(L;\mc{G}))$ with map from the latter to the former
induced by the coefficient homomorphism $\mc{IH}_{q}(L;\mc{G})\to
\mc{IH}_{q}(cL;c\mc{G})$.

For $I_{\pi}$, we define the collections
$\{\sigma^i_{\alpha}\}$, $\{\hat\sigma^i_{\alpha}\}$,
$\{\td\sigma^i_{\alpha}\}$, $\{\bd \hat \sigma^i_{\alpha}\}$, and $\{\bd
\td\sigma^i_{\alpha}\}$ as follows: Let $\{\sigma^i_{\alpha}\}$ denote the
collection of $i$-simplices of $B$, let $\hat\sigma^i_{\alpha}$ denote
the open neighborhood of $\sigma^i_{\alpha}$ in the second derived subdivision,
and let $\td\sigma^i_{\alpha}=\Pi^{-1}(\hat\sigma^i_{\alpha})$. For
each $\sigma^i_{\alpha}$, let $\bd \hat \sigma^i_{\alpha}$ denote $\hat
\sigma^i_{\alpha}\cap B^{i-1}$, i.e. the neighborhood of the boundary of
$\hat\sigma^i_{\alpha}$, and let $\bd \td\sigma^i_{\alpha}$ denote
$\td\sigma^i_{\alpha}\cap J^{i-1}=\Pi^{-1}(\bd \hat \sigma^i_{\alpha})$.
We claim that $H_{i}(IC_*(J^p,J^{p-1});\mc{G})\cong \oplus_{\alpha}
IH_i(\td\sigma^p_{\alpha},\bd
\td\sigma^p_{\alpha};\mc{G})$ for all $i$. The claim is obvious for $p=0$. For
$p>0$, we proceed by
induction on the simplices. In particular, order the simplices of each
dimension by placing an
ordering on the indices $\alpha$, and let $A^p_j$ equal the union of the
first $j$ of the
$\td\sigma^p_{\alpha}$. We will show by induction that
$IH_{i}(A^p_j,A^p_j\cap
J^{p-1};\mc{G})=\oplus_{\alpha\leq j} IH_i(\td\sigma^p_{\alpha},\bd
\td\sigma^p_{\alpha};\mc{G})$ for all
$j$, which will imply the claim since $A^p_j\subset A^p_{j+1}$ for all $j$,
$\cup A^p_j=J^p$, and
$J^p\cap J^{p-1}=J^{p-1}$. The formula obviously holds for $j=1$. For $j>1$,
consider the
Mayer-Vietoris sequence

\begin{multline*}
\to IH_i(A^p_j\cap \td\sigma^p_{j+1}, A^p_j\cap J^{p-1}\cap
\td\sigma^p_{j+1};\mc{G})\\ 
\to IH_{i}(A^p_j,A^p_j\cap J^{p-1};\mc{G})\oplus
IH_{i}(\td\sigma^p_{j+1},\bd
\td\sigma^p_{j+1};\mc{G})\\
\to IH_{i}(A^p_{j+1},A^p_{j+1}\cap
J^{p-1};\mc{G})\to.
\end{multline*}

\noindent Now, we  have $A^p_j\cap \td\sigma^p_{j+1}\subset
J^{p-1}$: In the
base $B$, the intersection
of $\cup_{i=0}^j \sigma^p_i$ and $\sigma^p_{j+1}$ must lie in a lower
dimensional skeleton and the
intersection of their neighborhoods will lie in the neighborhood of this
lower dimensional skeleton
hence in $B^{p-1}$, and we can then consider the inverse images of these
sets under the
projection. Therefore, the intersection terms of the sequence are $
IH_i(A^p_j\cap
\td\sigma^p_{j+1}, A^p_j\cap J^{p-1}\cap \td\sigma^p_{j+1};\mc{G})=0$. This
suffices to prove the
claim by a possibly infinite induction. Note, by the way, that it is in order
to employ such Mayer-Vietoris sequences that we must use the open sets $\hat
\sigma$ and their inverse images instead of simply the closed simplices, because Mayer-Vietoris sequences for intersection homology require open sets.

The corresponding results hold for $I_{\pi}-X_k\sim_{sphe}\hl_s(N,X_k)$ by
considering the analogously defined $ \bar
\sigma^p_{\alpha}=\pi^{-1}(\hat\sigma^p_{\alpha})$, etc., and in this case
we obtain $ H_{i}(IC_*(\bar J^p,\bar J^{p-1});\mc{G})\cong \oplus_{\alpha}
IH_i(\bar\sigma^p_{\alpha},\bd \bar\sigma^p_{\alpha};\mc{G})$.
Now we must examine the
$IH_i(\td\sigma^p_{\alpha},\bd \td\sigma^p_{\alpha};\mc{G})$ and
$IH_i(\bar\sigma^p_{\alpha},\bd \bar\sigma^p_{\alpha};\mc{G})$.

By definition, $\bar\sigma^p_{\alpha}=\pi^{-1}(\hat\sigma^p_{\alpha})$
and $\td\sigma^p_{\alpha}=\Pi^{-1}(\hat\sigma^p_{\alpha})$, which is the
mapping cylinder of the restriction of $\pi$ to $\bar\sigma^p_{\alpha}$.
According to the proof of \cite[Theorem 2]{Q2}, this restriction is precisely
the projection $\hl_s((N-B)\cup
\hat\sigma^p_{\alpha},
\hat\sigma^p_{\alpha})\to
\hat\sigma^p_{\alpha}$, which is a stratified fibration by
\cite[Cor. 6.2]{Hug}. By Corollary
\ref{C: strong triv}, above, there is a strong stratified trivialization
$\mf{h}_{\alpha}:\hat\sigma^p_{\alpha}\times \mc
L\to\hl_s(N,\hat\sigma^p_{\alpha})=\bar \sigma_{\alpha}^p$,
where $\mc L$ is the stratified fiber over some point of
$\hat\sigma^p_{\alpha}$. By  Corollary \ref{C: cyl strong triv}, we obtain a
strong stratum-preserving fiber homotopy equivalence
$h_{\alpha}:\hat\sigma^p_{\alpha}\times \bar
c\mc{L}\to\td\sigma^p_{\alpha}$. By
restricting the equivalences to the inverse images
$\pi^{-1}(\bd \hat\sigma^p_{\alpha})$, we see that the pairs
$(\bar\sigma^p_{\alpha},\bd \bar\sigma^p_{\alpha})$ and
$(\td\sigma^p_{\alpha},\bd \td\sigma^p_{\alpha})$ are stratum-preserving fiber homotopy equivalent to the pairs
$(\hat\sigma^p_{\alpha},\bd \hat\sigma^p_{\alpha})\times \mc{L}$ and
$(\hat\sigma^p_{\alpha},\bd \hat\sigma^p_{\alpha})\times\bar c \mc{L}$. Also, 
if $e_{\alpha}$ is the basepoint of $\hat\sigma^p_{\alpha}$ over which the the
trivializations are identities and $\mc{L}_{\alpha}$ is the corresponding fiber,
then, by bundle theory, $\mf{h}_{\alpha}^*\mc{G}\cong \hat
\sigma^p_{\alpha}\times (\mc{G}|\mc{L}_{\alpha})$ and 
$h_{\alpha}^*\mc{G}\cong \hat 
\sigma^p_{\alpha}\times (\mc{G}|\bar c\mc{L}_{\alpha})$.
Hence by the K\"{u}nneth theorem for the intersection homology of a product one of
whose factors is a manifold \cite{Ki} (and its obvious generalizations to
manifold pairs and a local coefficient system on the filtered space) and by the
homotopy equivalence $(\hat \sigma^p_{\alpha},\bd \hat
\sigma^p_{\alpha})\sim_{h.e.} (\sigma^p_{\alpha},\bd \sigma^p_{\alpha})$, we
obtain the formulas

{\small\begin{align*}
IH_i(\bar\sigma^p_{\alpha},\bd \bar\sigma^p_{\alpha};\mc{G})&\cong
IH_i((\hat\sigma^p_{\alpha},\bd \hat\sigma^p_{\alpha})\times
\mc{L}; h_{\alpha}^*\mc{G}) \cong
IH_i((\hat\sigma^p_{\alpha},\bd \hat\sigma^p_{\alpha})\times
\mc{L}; \hat
\sigma^p_{\alpha}\times (\mc{G}|\mc{L}_{\alpha}))\\
&\cong
H_p(\hat \sigma^p_{\alpha},\bd
\hat \sigma^p_{\alpha})\otimes IH_{i-p}(\mc{L} ;\mc{G}) \cong
H_p(\sigma^p_{\alpha},\bd
\sigma^p_{\alpha})\otimes IH_{i-p}(\mc{L};\mc{G})\\
IH_i(\td\sigma^p_{\alpha},\bd \td\sigma^p_{\alpha};\mc{G})&\cong
IH_i((\hat\sigma^p_{\alpha},\bd \hat\sigma^p_{\alpha})\times \bar c
\mc{L}; \mf{h}_{\alpha}^*\mc{G}) \cong
IH_i((\hat\sigma^p_{\alpha},\bd \hat\sigma^p_{\alpha})\times \bar c
\mc{L}; \hat
\sigma^p_{\alpha}\times (\mc{G}|\bar c\mc{L}_{\alpha}))\\
&\cong
H_p(\hat \sigma^p_{\alpha},\bd
\hat \sigma^p_{\alpha})\otimes IH_{i-p}(\bar c\mc{L} ;\mc{G})\cong
H_p(\sigma^p_{\alpha},\bd
\sigma^p_{\alpha})\otimes IH_{i-p}(\bar c\mc{L};\mc{G}).
\end{align*}}

\noindent (If $\mc{G}$ is a bundle of $R$ modules for some commutative ring $R$, then technically terms $H_p(\hat \sigma^p_{\alpha},\bd
\hat \sigma^p_{\alpha})$ should be take with $R$ coefficients, but we will suppress this from the notation for simplicity.) These are cellular simplices with coefficients in $ IH_{i-p}(\mc{L};\mc{G})$
and $ IH_{i-p}(\bar c\mc{L};\mc{G})$ (recall that on these spaces we are  using $\mc{G}$ to stand for $g^*\mc{G}$). We want to see that $IH_*(c\mc L;g^*\mc G)\cong IH_*(cL;\mc G)$ and $IH_*(\mc L;g^*\mc G)\cong IH_*(L;\mc G)$. 

We first show that $\bar c\mc{L}\sim_{s.p.h.e.}cL$ (recall that $L$ is our
geometric link). This will also imply $\mc{L}\sim_{s.p.h.e.}L$ since the
stratum preserving condition implies that $c\mc{L}-*\sim_{s.p.h.e.}cL-*$,
where $*$ represents the cone point, and clearly $\bar c\mc{L}-*\sim_{s.p.h.e.}c\mc{L}-*
\sim_{s.p.h.e.}\mc{L}$ and $cL-* \sim_{s.p.h.e.}L$. Once again we have
$\td\sigma^p_{\alpha}=\Pi^{-1}(\hat\sigma^p_{\alpha})$, which is the mapping
cylinder of the restriction of $\pi$ to $\bar\sigma^p_{\alpha}$, and this
restriction is precisely $\hl_s((N-X_k)\cup\hat\sigma^p_{\alpha},\hat\sigma^p_{\alpha})$.
Let us assume, as we are free to have imposed earlier, that we have chosen a
fine
enough triangulation so that each $\hat\sigma^i_{\alpha}$
has a
neighborhood in $N$ that is PL-homeomorphic to $\hat\sigma^i_{\alpha}\times cL$;
this is possible since each point in $X_k$ has a neighborhood of the form
$D^k\times cL$, and we can use generalized subdivisions (see e.g. \cite[\S
16]{MK}) to find a fine enough triangulation so that each
$\hat\sigma^i_{\alpha}$
lies in some such $D^k$. Then the restriction of $D^k\times cL$ to the
product over $\hat\sigma^i_{\alpha}$ provides such a neighborhood, which we
will denote $N^{\alpha}$.

Now, define a continuous function $\delta: \hat\sigma^p_{\alpha}\to
(0,\infty)$ such that
for each $x\in \hat\sigma^p_{\alpha}$, $B(x,\delta(x))\subset N^{\alpha}$,
where $B(x,\delta(x))$ is the closed ball in $N$ of radius $\delta(x)$
centered at $x$. For example, the function given by one half the distance from
the
boundary of the closure of $N^{\alpha}$ in $X$ will suffice. Then, following
\cite[p.
453]{Q1}, we can define
$\hl^{\delta}((N-X_k)\cup\hat\sigma^p_{\alpha},\hat\sigma^p_{\alpha})$ to be
the subset of
$j\in\hl((N-X_k)\cup\hat\sigma^p_{\alpha},\hat\sigma^p_{\alpha})$ such that
Im$(j:I\to N)$ lies within $\delta(j(0))$ of $j(0)$. Similarly, we can define
$\hl_s^{\delta}((N-X_k)\cup\hat\sigma^p_{\alpha},\hat\sigma^p_{\alpha})$.
According to \cite[Lemma 2.4]{Q1}, the inclusion
$\hl^{\delta}((N-X_k)\cup\hat\sigma^p_{\alpha},\hat\sigma^p_{\alpha})\into
\hl((N-X_k)\cup\hat\sigma^p_{\alpha},\hat\sigma^p_{\alpha})$ is a fiber
homotopy equivalence. Since the proof of this lemma consists of appropriately ``shrinking''
each path of the holink along itself, restriction provides a stratum-preserving homotopy
equivalence of
$\hl_s^{\delta}((N-X_k)\cup\hat\sigma^p_{\alpha},\hat\sigma^p_{\alpha})$ and
$\hl_s((N-X_k)\cup\hat\sigma^p_{\alpha},\hat\sigma^p_{\alpha})$. On the other
hand, by our choice of $\delta$, clearly
$\hl_s^{\delta}((N-X_k)\cup\hat\sigma^p_{\alpha},\hat\sigma^p_{\alpha})=
\hl_s^{\delta}(N^{\alpha},\hat\sigma^p_{\alpha})$, and by the same arguments, the inclusion
$\hl_s^{\delta}(N^{\alpha},\hat\sigma^p_{\alpha})\into \hl_s(N^{\alpha},\hat\sigma^p_{\alpha})$ must be a stratum-preserving fiber
homotopy equivalence. Hence, we have shown that
$\hl_s((N-X_k)\cup\hat\sigma^p_{\alpha},\hat\sigma^p_{\alpha})$ is stratum-preserving fiber homotopy
equivalent to $\hl_s(N^{\alpha},\hat\sigma^p_{\alpha})$, and this induces a
stratum-preserving fiber homotopy equivalence of the mapping cylinders of
each of the projections to $\hat\sigma^p_{\alpha}$ by Corollary \ref{C: cyl
fib 
hom eq}. By definition, we already
know that the mapping cylinder of the projection of $\hl_s((N-X_k)\cup\hat\sigma^p_{\alpha},\hat\sigma^p_{\alpha})$ is $\td \sigma_{\alpha}$, and by the same methods as our
original proof that $I_{\pi}\sim_{s.p.h.e.}N$, we can show that $N^{\alpha}$
is stratum-preserving homotopy equivalent (by a map we will denote $\phi$) to the mapping cylinder of the projection 
$\pi_{N^{\alpha}}:\hl_s(N^{\alpha},\hat\sigma^p_{\alpha})\to
\hat \sigma^p_{\alpha}$ (as our nearly stratum-preserving deformation
retraction
of $N^{\alpha}\cong \hat \sigma^p_{\alpha}\times cL$ to $ \hat
\sigma^p_{\alpha}$, we can take the product with the identity on  $ \hat 
\sigma^p_{\alpha}$ of the obvious nearly stratum-preserving retraction of $cL$
to
it cone point). Hence
we have $N^{\alpha}\sim_{s.p.h.e}\td\sigma^p_{\alpha}$. Finally, we know that $\td\sigma^p_{\alpha}\sim_{s.p.h.e.}
\hat\sigma^p_{\alpha}\times \bar c\mc{L}$, while
$N^{\alpha}\cong_{P.L}\hat\sigma^p_{\alpha}\times cL$ by a stratum-preserving
homeomorphism. Since $\hat\sigma^p_{\alpha}$ is contractible and unfiltered, we at last
conclude:  
\begin{equation*} 
\bar c\mc{L}\sim_{s.p.h.e.} \hat\sigma^p_{\alpha}\times \bar c\mc{L}
\sim_{s.p.h.e.} \td\sigma^p_{\alpha} \sim_{s.p.h.e.}
N^{\alpha}\cong_{P.L.}\hat\sigma^p_{\alpha}\times cL \sim_{s.p.h.e.} cL.  
\end{equation*}

In order to show that $IH_*(c\mc L;g^*\mc G)\cong IH_*(cL;\mc G)$, we must
recall that the bundle over $\mc L$ is actually that determined by the
pullback of $\mc{G}$ under the stratum-preserving homotopy equivalence
$g:I_{\pi}\to N$ and
show that the homotopy equivalences of the preceding paragraph are covered by
the proper coefficient maps. In particular, for any $\alpha$, let us
identify $cL$ with one of its isomorphic images in $N^{\alpha}$ under the
inclusion $*\times cL\into \hat \sigma_{\alpha}\times cL$, for some $*$,
followed by the isomorphism $\hat \sigma_{\alpha}\times cL\to N^{\alpha}$.
Let $h:N^{\alpha}\to \td \sigma_{\alpha}$ denote the stratum-preserving
homotopy
equivalence
determined as in the preceding paragraph. Then $h^*g^*\mc{G}=(gh)^*\mc{G}$ is a
coefficient
system on $N^{\alpha}$. We want to see that it is equivalent to
$\mc{G}|N^{\alpha}$. 
Let $\pi_{N^{\alpha}}:\hl_s(N^{\alpha}, \hat \sigma^p_{\alpha})\to  \hat
\sigma^p_{\alpha}$, $\pi^{\delta}_{N^{\alpha}}:\hl^{\delta}_s(N^{\alpha},
\hat
\sigma^p_{\alpha})\to \hat \sigma^p_{\alpha}$, and
$\pi^{\delta}:\hl_s^{\delta}((N-X_k)\cup\hat\sigma^p_{\alpha},\hat\sigma^p_{\alpha})\to \hat \sigma^p_{\alpha}$ be the projections. Let
$\psi: \hl_s(N^{\alpha}, \hat \sigma^p_{\alpha})\to \hl^{\delta}_s(N^{\alpha},
\hat
\sigma^p_{\alpha})$ be the stratum-preserving fiber homotopy inverse to the
inclusion. By the proof of \cite[Lemma 2.4]{Q1}, if $\gamma\in  \hl_s(N^{\alpha},
\hat \sigma^p_{\alpha})$, then $\psi\gamma:[0,1]\to N^{\alpha}$ is the path
$\gamma$ composed with multiplication by a number $r\in (0,1)$, which may depend on
$\gamma$. Then $\psi$ induces a map $\bar\psi: I_{\pi_{N^{\alpha}}}\to
I_{\pi^{\delta}_{N^{\alpha}}}$. Also, the inclusion $\hl^{\delta}_s(N^{\alpha},
\hat
\sigma^p_{\alpha})\into \hl_s((N-X_k)\cup\hat\sigma^p_{\alpha},\hat\sigma^p_{\alpha})$ induces an inclusion $I_{\pi^{\delta}}\to I_\pi$, and the
equality $\hl_s^{\delta}(N^{\alpha},\hat
\sigma^p_{\alpha})=\hl_s^{\delta}((N-X_k)\cup\hat\sigma^p_{\alpha},\hat\sigma^p_{\alpha})$
induces an equality $I_{\pi^{\delta}_{N^{\alpha}}}=I_{\pi^{\delta}}$. With this
notation, $gh$ is the composition of maps
\begin{equation*}
\begin{CD}
N^{\alpha}&@>\phi >>  I_{\pi_{N^{\alpha}}}&@>\bar \psi>>&
I_{\pi^{\delta}_{N^{\alpha}}}=I_{\pi^{\delta}}\into \td
\sigma^p_{\alpha}\into I_{\pi}&@>g>>& N
\end{CD}
\end{equation*}
(recall that  $\td 
\sigma^p_{\alpha}$ is the mapping cylinder of the projection $\pi|_{\pi^{-1}(
\hat \sigma^p_{\alpha})}=\pi|_{\hl_s((N-X_k)\cup\hat\sigma^p_{\alpha},\hat\sigma^p_{\alpha})}$ ).
If $d(\cdot,\cdot)$ represents distance in $X$ and $(s,t,x)\in \hat \sigma^p_{\alpha}\times (0,1)\times L\subset \hat
\sigma^p_{\alpha}\times cL-\hat \sigma^p_{\alpha}$ and $\gamma_{(s,t,x)}$ is
the
path of the retraction of $(s,t,x)$ to its cone point, $s\in\hat
\sigma^p_{\alpha}$, under the nearly
stratum-preserving deformation retraction of $N^{\alpha}$ to $\hat
\sigma^p_{\alpha}$, then we can calculate $gh(s,t,x)$ as follows:
$\phi(s,t,x)=(d((s,t,x),X_k), \gamma_{(s,t,x)})\in I_{\pi_{N^{\alpha}}}$; $\bar \psi
\phi (s,t,x)=(d((s,t,x),X_k), \psi\gamma_{(s,t,x)})$; the next three maps are
inclusions; and finally $gh(s,t,x)= S(\psi\gamma_{(s,t,x)},d((s,t,x),X_k))(d((s,t,x),X_k))$,
which
lies on the line from $(s,t,x)$ to the cone point $s\in \hat
\sigma^p_{\alpha}$ (recall that we have assumed that all point in $N$ are distance $<1$ from $X_k$ and $S$ retracts all paths along themselves). Therefore, if, for a given $t\in (0,1)$, $H:\hat
\sigma^p_{\alpha}\times (0,1)\times
L\times I\to \hat \sigma^p_{\alpha}\times (0,1)\times
L$ is any of the  stratum-preserving deformation retractions to $\hat
\sigma^p_{\alpha}\times
t\times L$  which takes each $s\times (0,1)\times x$ to
$(s,t,x)$, then $H\circ(gh|_{\hat
\sigma^p_{\alpha}\times t\times L}\times \text{id}_I)$ is a homotopy from
$gh|_{\hat
\sigma^p_{\alpha}\times t\times L}$ to the inclusion $i:\hat
\sigma^p_{\alpha}\times t\times L\into N^{\alpha}$. Therefore, $(gh|_{\hat
\sigma^p_{\alpha}\times t\times L})^*\mc{G}\cong \mc{G}|_{\hat
\sigma^p_{\alpha}\times t\times L}$. But then, employing some bundle
theory
and recalling that $\mc{G}$ is only defined over the top stratum, we have
\begin{equation*}
\mc{G}\cong \mc{G}|_{\hat \sigma^p_{\alpha}\times t\times L}\times (0,1)\cong
(gh|_{\hat
\sigma^p_{\alpha}\times t\times L})^*\mc{G}\times (0,1)\cong
(gh)^*\mc{G}|_{\hat
\sigma^p_{\alpha}\times t\times L}\times (0,1)\cong (gh)^*\mc{G}.
\end{equation*}

Therefore, $IH_*^{\bar
p}(N^{\alpha};h^*g^*\mc{G})\cong IH_*^{\bar
p}(N^{\alpha};\mc{G}|N^{\alpha})$, and, by Corollary \ref{C: s.p.h.e. with l.c.},
we have isomorphisms $IH_*^{\bar p}(\td
\sigma_{\alpha};g^*\mc{G})\cong IH_*^{\bar p}(N^{\alpha};h^*g^*\mc{G})$.
Also, there are isomorphisms $IH_*^{\bar
p}(N^{\alpha};\mc{G}|N^{\alpha})\cong
IH_*^{\bar p}(cL;\mc{G}|cL)$ since the stratum-preserving homotopy equivalence
of $cL$ and
$N^{\alpha}$ is induced by inclusion. Similarly, we obtain $IH_*^{\bar p}(\td
\sigma_{\alpha};g^*\mc{G})\cong IH_*^{\bar p}(\bar c\mc{L};g^*\mc{G}|\bar 
c\mc{L})$ due
to the strong stratum-preserving trivialization $\hat\sigma_{\alpha}\times
c\mc{L}\to \td \sigma_{\alpha}$, which is the
identity on $*\times \mc{L}$, and the stratum-preserving homotopy equivalence
given by the
inclusion $\bar c\mc{L}\to \hat\sigma_{\alpha}\times \bar c\mc{L}$. 
Lastly, since $IH_*(c\mc{L};g^*\mc{G}) \cong IH_*(\bar{c} \mc{L};g^*\mc{G})$, induced by inclusion, we conclude  $IH_*(c\mc L;g^*\mc G)\cong IH_*(cL;\mc G)$.

In the same manner, we
can use
these homotopy equivalences to show that $IH_*(\mc L;g^*\mc G)\cong IH_*(L;\mc
G)$: the restriction $h|N^{\alpha}-\hat \sigma_{\alpha}$ is a stratum-preserving homotopy
equivalence to $\td \sigma_{\alpha}-\hat \sigma_{\alpha}$, and
$gh|N^{\alpha}-\hat \sigma_{\alpha}$ is stratum-preserving homotopic to the
inclusion. As above, this gives an isomorphism $ IH_*^{\bar p}(cL-*;\mc{G}|cL-*)
\cong IH_*^{\bar p}(\bar
c\mc{L}-*;g^*\mc{G}|\bar c\mc{L}-*)$. But $IH_*^{\bar p}(cL-*;\mc{G}|cL-*)\cong
IH_*^{\bar p}(L;\mc{G}|L)$ and $ IH_*^{\bar p}(\bar
c\mc{L}-*;g^*\mc{G}|\bar c\mc{L}-*)\cong IH_*^{\bar p}(\mc{L};g^*\mc{G}|\mc{L})$,
as induced by the inclusions (for the former, the inclusion can be considered as
$L\into L\times t\subset cL$ for any $t\in (0,1)$; also note that the
$\mc{G}|(L\times t)$ are all isomorphic, as $\mc{G}|cL\cong (\mc{G}|L\times
s)\times (0,1)$ for any $s\in (0,1)$). We should also observe that, as seen above
in our study of local coefficient systems on stratified fibrations (and their
mapping cylinders), these intersection homology groups are independent, up to
isomorphism, of our particular choice of fiber $\mc{L}$. 

At this point, by putting together the results of the preceding paragraphs, we
have shown that $E^1$ and $\bar E^1$ are composed of the correct cellular chain
modules. It remains to verify that the boundary maps are those compatible with the
cellular chain complex of $X_k$ with the desired local coefficient
systems and
that the map between the two systems is the desired one. This will be done over
the
following pages. 
\end{proof}

Let us first dispense with the lemma whose proof was deferred.
\begin{proof}[Proof of lemma \ref{l: def ret}]
According to \cite[Prop. 3.2]{Q1} (also see \cite[p. 241]{Q2}) and the
remarks preceding Definition A1 of \cite{Hug}), there exists a
neighborhood $V$ of $X_k$ in $N$ and a nearly stratum-preserving deformation retraction $r: V\times
I\to X$ of $V$ to $X_k$. Note that $r$ may take points in $V$ outside of $V$ during the
retraction process, which is the reason this lemma is necessary. Let $N'$ be another open regular neighborhood  of $X_k$ such that
$\bar N'\subset \text{int}(V)$. Such an $N'$ exists because we can find a closed
regular neighborhood  of $X_k$ in int$(V)$ (as by the construction of \cite[Prop. 1.11]{Sto}), but this
will clearly also be a regular neighborhood of $X_k$ in $X$. Furthermore, we have
$\bar N'\subset N$. By \cite[Prop. 1.5]{Sto}, cl$(\bar N-\bar N')$ is PL-isomorphic to the
product fr$(\bar N')\times I$, and by an isomorphism which takes each cl$(X_i\cap
(\bar N-\bar N'))$ to $(X_i\cap \text{fr}(\bar N'))\times I$. Hence, there is a
homotopy which gives a deformation retraction of $N$ to $\bar N'$ which is
stratum preserving and maps $N$ into $N$, namely the homotopy which retracts cl$(\bar
N-\bar N')\cong  \text{fr}(\bar N')\times I$ to $ \text{fr}(\bar N')\times \{0\}$ while leaving the
rest of $N$ fixed. So, if we consider $N$ as $N'\cup_{\text{fr}(\bar N')\times \{0\}}\text{fr}(\bar N')\times I$ and represent each point in $\text{fr}(\bar N')\times I$ by $(y,s)$, we define
\begin{equation*}
H(t,x)=
\begin{cases}
x,& x\in N- \text{cl$(\bar N-\bar N')$}\\
(y,ts), & x\in\text{cl}(\bar N-\bar N'), x=(y,s). 
\end{cases}
\end{equation*}
This is clearly continuous.
Now we can define $R: N\times I \to N$ as
\begin{equation*}
R(t,x)=
\begin{cases}
H(3t-2,x), & t\in[2/3,1]\\
H(0,x)=r(1, H(0,x)), & t\in[1/3,2/3]\\
r(3t,H(0,x)), & t\in[0,1/3].
\end{cases}
\end{equation*}
This is the homotopy which shrinks $N$ to $N'\subset V$ by deformation
retracting
$\text{fr}(\bar N')\times I$ as time decreases from $1$ to $2/3$, then
pauses in the time interval $[1/3,2/3]$, and then finally, for time in
$[0,1/3]$, deformation retracts $H(0,N')$ to $X_k$ utilizing the retraction
$r$
on $V$. Furthermore, $R$ is nearly stratum-preserving since $H$ is stratum-preserving and $r$ is nearly stratum-preserving, and continuity follows from the continuity of $H$ and $r$. Hence $R$
is our desired deformation retraction.  
\end{proof}

We will now complete the proof of Proposition \ref{P: Neigh IH} by showing that the
$E^1$ terms constructed above together with the boundary maps given by the spectral sequences are chain isomorphic to cellular chain complexes on $B$ with local coefficients determined as above from the stratified fibration. The proof is essentially that given in the usual computation of the $E^2$ term of the spectral sequence of a fibration (see, e.g., Whitehead \cite[Chapter XIII]{Wh}). We run through the details in this situation mainly to double check that our departures are consistent with the usual the proof. We make the following initial remarks with the hope of clarifying the specific points of departure:
\begin{enumerate}
\item We are using intersection homology as opposed to the more traditional homology theories of the standard presentations. This technicality is dealt with by the use
of stratified fibrations. We have already employed the stratified fibration
$\hl_s(X,X_k)\to X_k$ to construct both bundles of coefficients. However, in the
following, as in
the standard proof, it is the local triviality that comes to the fore.

\item We need to allow a bundle of local coefficients on the total space $N$
(or, more specifically, on the complement of $N\cap X_{n-2}$). For ordinary homology
this
complication is tacitly present in Whitehead's treatment for extraordinary
homology theories. Similarly, local coefficients fit more-or-less seamlessly into
the following, only warranting an occasional remark to keep us honest.

\item Due to our need for all the subsets to be open in order to employ the
Mayer-Vietoris arguments for intersection homology used in the first part of the
proof, we have
been saddled with using the groups $H_*(\hat \sigma_{\alpha},\bd\hat
\sigma_{\alpha})$ as the generators of our cellular chain complex on $X_k$.
However, this presents no difficulty as the homotopy equivalences $(\hat
\sigma_{\alpha},\bd\hat \sigma_{\alpha})\sim_{h.e.}(\sigma_{\alpha},\bd
\sigma_{\alpha})$ induce canonical chain isomorphisms between the cellular chain
complex so-defined and that of the standard cellular chain complex. After a brief
initial statement, below, we will usually suppress the distinction.

\item In contrast to the above generality, the following discussion is simplified by the
fact that we use a fixed triangulation of $X_k$ to impart its structure as a
CW-complex. Hence, we can treat all characteristic maps as inclusions (although, once again, we take the unusual approach of including the open simplices $(\hat\sigma_{\alpha},\bd\hat\sigma_{\alpha})$).
\end{enumerate}

\begin{proof}[Completion of the proof of Theorem \ref{P: Neigh IH}]
We will employ the following conventions: For the homology of $X_k=B$ with local
coefficients, we will employ the paradigm of computing homology with local
coefficients by using coefficients in the module over
a basepoint of each relative singular simplex generating each cellular simplex
(as
opposed to defining the coefficient as a lift as we did above). In this manner, we
will define cellular homology of $B$ with local coefficients as in \cite{Wh}. We
use the given triangulation of $B$ to define it as a (regular) cell complex.
Hence, with local coefficient system $\mc{M}$ with fiber $M$, we define
$H_*(B;\mc{M})$ as the
homology of the chain complex of $B$ with local coefficients
\begin{align*}
\Gamma_*(B;\mc{M})&=
\oplus_{\alpha}H_*(\hat\sigma_{\alpha},\bd \hat\sigma_{\alpha};\mc{M}) \cong
H_*(B^*,B^{*-1};\mc{M})\\
&\cong
\oplus_{\alpha}H_*(\sigma_{\alpha},\bd \sigma_{\alpha};\mc{M}) 
\cong \oplus_{\alpha}M,
\end{align*}
with the boundary maps given by those of the the long exact sequence of
the
triple $(B^*, B^{*-1},
B^{*-2})$ with coefficient maps handled as usual for homology with local
coefficients (see \cite{Wh}). Again, it is the  direct sum
$\oplus_{\alpha}H_*(\sigma_{\alpha},\bd \sigma_{\alpha};\mc{M})$
which is
generally used to define
the cellular homology chain complexes, but the above string of isomorphisms and
those induced in the obvious places by the homotopy equivalences of the
$\sigma_{\alpha}$ and the $\hat\sigma_{\alpha}$ allow us to redefine the cellular
homology in this manner. Note that each $H_*(\hat\sigma_{\alpha},\bd
\hat\sigma_{\alpha};\mc{M})\cong H_*(\sigma_{\alpha},\bd \sigma_{\alpha};\mc{M})$
can be considered as generated by the singular simplex
given by the characteristic inclusion $i_{\alpha}:\Delta^i\to
\sigma_{\alpha}$ with coefficient in the module over the leading vertex
$e_{\alpha}:=i_{\alpha}(e_0)$.

We have already seen that there is an isomorphism 
\begin{equation*}
\Gamma_p(B;\mc{IH}^{\bar p}_q(\bar c \mc L;\mc{G}|\bar c \mc
L)) \overset{\cong}{\to} IH_{p+q}^{\bar p}(J^p, J^{p-1};\mc{G}|J^p),
\end{equation*} 
in which $\bar c \mc L$ is the closed cone on the fiber of the stratified
fibration
$\pi:\hl_s(N,X_k)\to B=X_k$ up to stratum-preserving homotopy equivalence.  Let us
specify this isomorphism, which we will denote $\lambda$, by determining its
action on the cellular simplex with coefficient $x[i_{\alpha}]$, where $[i_{\alpha}]$
denotes the
homology class of the generating singular chain $i_{\alpha}$ in $H_p(\hat
\sigma_{\alpha}, \bd \hat \sigma_{\alpha})$ and $x\in IH_q^{\bar p}(\bar c \mc
L_{\alpha};\mc{G})$, $ \bar c \mc L_{\alpha}$ the cone on the fiber over
$e_{\alpha}=i_{\alpha}(e_0)$. Let $h_{\alpha}$ be the strong stratified
trivialization $h_{\alpha}:(\hat\sigma_{\alpha}, \bd \hat \sigma_{\alpha})\times
\bar c \mc L_{\alpha} \to (\td\sigma_{\alpha}, \bd\td \sigma_{\alpha})$ obtained
in the first part of the proof from a strong stratified
trivialization from $\hat \sigma_{\alpha}\times \mc{L}_{\alpha}$ to
$\hl_s(N,\hat \sigma_{\alpha})$ (and covered by the induced map
$h_{\alpha}^*(\mc{G}|\td\sigma_{\alpha})\to(\mc{G}|\td\sigma_{\alpha})$); and let
$[i_{\alpha}]\times :IH_q^{\bar p}( \bar c \mc L_{\alpha};\mc{G})\to
IH_{p+q}^{\bar p}(\hat
\sigma_{\alpha}\times \bar c \mc L, \bd\hat\sigma_{\alpha}\times \bar c \mc
L;h_{\alpha}^*\mc{G})$ be the map given by the cross product of the K\"{u}nneth
theorem and using the unique isomorphism $\hat \sigma_{\alpha}\times
(\mc{G}|\bar c\mc{L}_{\alpha})\overset{\cong}{\to}h_{\alpha}^*\mc{G}$ which is the identity over $e_{\alpha}$. Let
$j_{\alpha}:(\td \sigma_{\alpha}, \bd \td \sigma_{\alpha})\to (J^p, J^{p-1})$ be
the inclusion. Then, by reversing the order of the steps by which we previously
demonstrated this isomorphism, we see that $\lambda$ can be given by 
\begin{equation*}
\lambda(x[i_{\alpha}])=j_{\alpha*}h_{\alpha*}([i_{\alpha}]\times 
x). 
\end{equation*}

To show that, for each fixed $q$, $E^1_{*,q}$ with boundary $d^1$ is chain
isomorphic to the complex
$\Gamma_p(B;\mc{IH}^{\bar p}_q(\bar c \mc L;\mc{G}|\bar c \mc
L))$, it remains to prove that the following diagram commutes for each $p$:
\begin{diagram}\label{D: lambda com}
\Gamma_p(B;\mc{IH}^{\bar p}_q(\bar c \mc L;\mc{G}|\bar c \mc L)) & \rTo^{\bd}
&\Gamma_{p-1}(B;\mc{IH}^{\bar p}_q(\bar c \mc L;\mc{G}|\bar c \mc L))\\
\dTo^{\lambda} & & \dTo^{\lambda}\\
IH_{p+q}^{\bar p}(J^p, J^{p-1};\mc{G}|J^p) & \rTo^{\bd} & IH_{p+q-1}^{\bar
p}(J^{p-1}, J^{p-2};\mc{G}|J^p).\\
\end{diagram}
(Note: we use the symbol ``$\bd$'' to refer to either ``boundary'' map.) Let
$\Gamma_*(\hat\sigma^p_{\alpha})$ be the cellular chain complex generate by
$H_*(\hat\sigma^p_{\alpha},\bd\hat\sigma^p_{\alpha})\cong
H_p(\hat\sigma^p_{\alpha},\bd\hat\sigma^p_{\alpha})$. Let $\Gamma_*(\bd
\hat\sigma^p_{\alpha})$ be the cellular chain complex induced on $\bd
\hat\sigma^p_{\alpha}$; in other words $\oplus_{\beta<\alpha}H_{p-1}(\hat
\sigma^{p-1}_{\beta},\bd\hat \sigma^{p-1}_{\beta})$, the sum over the 
facets of $\sigma_{\alpha}$; similarly for local coefficients. Let $i_{\alpha}:
\hat\sigma_{\alpha}\into B$ and $j_{\alpha}:  (\td \sigma^p_{\alpha},\bd \td
\sigma^p_{\alpha})\into (J^p, J^{p-1})$ be the inclusions. Let $\mc{L}_{\beta}$
and $\mc{L}_{\alpha}$ stand for the fibers over the basepoints of the given
generators of $H_{p-1}(\hat\sigma_{\beta},\bd\hat\sigma_{\beta})$ and
$H_{p}(\hat\sigma_{\alpha},\bd\hat\sigma_{\alpha})$. We consider the following
cube diagram whose front face is diagram \eqref{D: lambda com}:

\pagebreak
{\scriptsize
\begin{diagram}[PostScript=dvips,landscape,LaTeXeqno,balance]
\label{E: cube}
&\\
&\\
&\\
&\\
&\\
&\\
&\\
&\\
&\\
&\\
&\\
&\\
&\\
&\\
&\\
&\\
&\\
&\\
&\\
&\\
&\\
&\\
&\\
&\\
&\\
&\\
&\\
&\\
&\\
&\\
&\\
&\\
&\\
&\\
&\\
&\\
H_p(\hat\sigma^p_{\alpha}, \bd \hat\sigma^p_{\alpha};IH^{\bar p}_q(\bar
c
\mc L_{\alpha};\mc{G}|\bar c \mc L_{\alpha}) )    &   &
\rTo^{\bd}  &   & \oplus_{\beta<\alpha} H_{p-1}(\hat\sigma^{p-1}_{\beta},
\bd \hat\sigma^{p-1}_{\beta};IH^{\bar p}_q(\bar c \mc
L_{\alpha};\mc{G}|\bar c \mc L_{\alpha}) )   & &\\
& \rdTo_{i_{\alpha*}} & & & \vLine^{\lambda'} & \rdTo_{i_{\alpha*}} &\\
\dTo^{\lambda} & & \Gamma_p(B;\mc{IH}^{\bar p}_q(\bar c \mc L;\mc{G}|\bar c 
\mc L)) &\rTo^{\bd} & \HonV & &
\Gamma_{p-1}(B;\mc{IH}^{\bar p}_q(\bar c \mc L;\mc{G}|\bar c \mc L))\\
&&\dTo^{\lambda} & & \dTo & &\\
IH_{p+q}^{\bar p}(\td \sigma^p_{\alpha}, \bd\td \sigma^p_{\alpha};\mc{G}) & 
\hLine & \VonH & \rTo^{\bd} & IH_{p+q-1}^{\bar p}(\bd\td
\sigma^p_{\alpha},\bd\td \sigma^p_{\alpha}\cap
J^{p-2};\mc{G}) &&\dTo_{\lambda}\\
&\rdTo_{j_{\alpha*}} &&&&\rdTo_{j_{\alpha*}} &\\
&& IH_{p+q}^{\bar p}(J^p, J^{p-1};\mc{G}|J^p)&&\rTo^{\bd} &&
IH_{p+q-1}^{\bar p}(J^{p-1}, J^{p-2};\mc{G}|J^p)\\
\end{diagram}}
\pagebreak

The top
back map of the cube \eqref{E: cube} is the
boundary map with constant coefficients. The map
$\lambda'$ is
determined by the restriction of the strong trivialization over $\hat
\sigma^p_{\alpha}$, not by the trivialization given over each $\hat
\sigma_{\beta}$, i.e. $\lambda'(x[i_{\beta}])=\dot h_{\alpha *}([i_{\beta}]\times
x)$, where $\dot h_{\alpha}$ is 
the restriction of $h_{\alpha}$ to $\bd \hat \sigma_{\alpha}\times \bar c 
\mc{L}_{\alpha}$ and $[i_{\beta}]\times$ is the cross product of the K\"{u}nneth
theorem together with the bundle isomorphism  $[\hat \sigma_{\alpha}\times
(\mc{G}|\bar c\mc{L}_{\alpha})]|(\hat\sigma_{\beta}\times \mc{L}_{\alpha}) 
\overset{\cong}{\to}h_{\alpha}^*\mc{G}|(\hat\sigma_{\beta}\times
\mc{L}_{\alpha})$.
Note also that the $i_{\alpha*}$ on the top right edge is
induced by the inclusion on the spaces $\hat \sigma_{\beta}$ and the map on
coefficients determined by the map of the coefficient system over the path class 
in
$\hat \sigma_{\alpha}$
between the basepoints of the generators of
$H_p(\hat\sigma_{\alpha},\bd\hat\sigma_{\alpha})$ and
$H_{p-1}(\hat\sigma_{\beta},\bd\hat\sigma_{\beta})$.

The top and bottom of this cube are clearly commutative.  On the left side, it is
easy to see then that the difference
between the maps $j_{\alpha*}\lambda$ and $\lambda i_{\alpha*}$ is the difference
between first including into a direct sum and then mapping the direct sum
componentwise versus first mapping components and then including into the direct
sum. These are equivalent and yield commutativity of the left side of the
cube. The commutativity of the right side of the cube is harder to prove and will
be deferred until we complete the rest of the proof of Proposition \ref{P: Neigh
IH}.

Finally, we show that the back face of the cube commutes:
If $i_{\beta}$ are the facets of $i_{\alpha}$, suitably ordered, then by standard
cellular theory, $\bd [i_{\alpha}]=\sum_{\beta=0}^p
(-1)^{\beta}[i_{\beta}]$. Within $\td
\sigma_{\alpha}$, we have $\lambda(x[i_{\alpha}])= h_{\alpha*}([i_{\alpha}]\times
x)$. Now we compute:
\begin{align*}
\lambda'\bd(x[i_{\alpha}])&=\lambda'\sum_{\beta=0}^p (-1)^{\beta}x[i_{\beta}]\\
&= \sum_{\beta=0}^p (-1)^{\beta}\dot h_{\alpha*}([i_{\beta}]\times x)\\
&=\dot h_{\alpha*}\left(\sum_{\beta=0}^p (-1)^{\beta}[i_{\beta}]\times x\right)\\
&=\dot h_{\alpha*} \bd([i_{\alpha}]\times x)\\
&=\bd h_{\alpha*} ([i_{\alpha}]\times x)\\
&=\bd \lambda(x[i_{\alpha}]).
\end{align*}
The fourth equality is due to the naturality of the K\"{u}nneth theorem and, more particularly, the Eilenberg-Zilber map for intersection homology (\cite{Ki}). 

We have shown the commutativity of all faces of the cube but for the front face
(and the deferred right face, which we now assume to commute). Thus we have
\begin{align*}
\lambda\bd i_{\alpha*}&=\lambda(i_{\alpha}|\bd \hat \sigma_{\alpha})_* \bd\\
&=(j_{\alpha}|\bd \td \sigma_{\alpha})_*\lambda'\bd\\
&=(j_{\alpha}|\bd \td \sigma_{\alpha})_*\bd\lambda\\
&=\bd j_{\alpha*}\lambda\\
&=\bd\lambda i_{\alpha*}.
\end{align*} 
Hence, we have established commutativity of diagram \eqref{D: lambda com} on every term of the direct sum
decomposition of
$\Gamma_p(B; \mc{IH}_q^{\bar p}(\bar c\mc L;\mc G|\bar c\mc L))$ and therefore
commutativity in general. So we have shown that $E^1_{*,q}$ and
$\text{Cell}_*(X_k;\mc{IH}_q^{\bar p}(cL;\mc G))$ are chain equivalent as desired
(recalling once again that $IH_*^{\bar p}(L;\mc{G})\cong IH_*^{\bar
p}(\mc{L};\mc{G})$). The proof for $\bar
E^1_{*,q}$ and $\text{Cell}_*(X_k;\mc{IH}_q^{\bar p}(L;\mc G))$ is entirely
analogous using the $\bar \sigma_{\alpha}$ instead of the $\td \sigma_{\alpha}$.

Finally, we need to see that the map $i_*:\bar E^1_{*,*}\to E^1_{*,*}$
which is induced by
the inclusion $i: \hl_s(N,X_k)\into I_{\pi}$ is given by the cellular chain map
$j: \text{Cell}(X_k,\mc{IH}_q^{\bar p}(\mc L;\mc G))\to$ Cell$(X_k,\mc{IH}_q^{\bar
p}(\bar c\mc L;\mc G))$ induced by the map on the bundle of coefficients that is
induced on the fiber intersection homology by the inclusion $\mc
L\into \bar c\mc L$. In other words, we must show that the following diagram
commutes:
\begin{equation*}
\begin{CD}
\Gamma_*(B;\mc{IH}^{\bar p}_q( \mc L;\mc{G}| \mc L)) &@>j_*>> &\Gamma_*(B;\mc{IH}^{\bar p}_q(\bar c \mc L;\mc{G}|\bar c \mc L)) \\
@V\lambda VV&&@V\lambda VV\\
IH_{*+q}^{\bar p}(\bar J^*, \bar J^{*-1};\mc{G}|\bar J^*) & @>i_* >> & IH_{*+q}^{\bar
p}(J^*, J^{*-1};\mc{G}|J^*).
\end{CD}
\end{equation*}
But this commutes for each $*$ as follows directly by considering how the maps $\lambda$
act on each summand of $\Gamma_*(B;\mc{IH}^{\bar p}_q( \mc L;\mc{G}| \mc L))$ or
$\Gamma_*(B;\mc{IH}^{\bar p}_q(\bar c \mc L;\mc{G}|\bar c \mc L))$ determined by the
cell $\hat \sigma_{\alpha}$. It is only necessary to employ the naturality of the
K\"{u}nneth theorem and the fact that each trivializing stratum-preserving fiber homotopy equivalence
$h_{\alpha}:(\hat\sigma_{\alpha}, \bd \hat \sigma_{\alpha})\times \bar c \mc L \to
(\td\sigma_{\alpha}, \bd\td \sigma_{\alpha})$ as defined restricts to a stratum-preserving fiber homotopy equivalence $\mf
h_{\alpha}:(\hat\sigma_{\alpha}, \bd \hat \sigma_{\alpha})\times \mc L \to
(\bar\sigma_{\alpha}, \bd\bar \sigma_{\alpha})$ (in fact $h_{\alpha}$ was defined by
extending $\mf h_{\alpha})$.

We have already seen that $\lambda$ is a chain map. That $i_*$ and $j_*$ are, in fact, maps of chain complexes follows from the naturality of homology and intersection homology with respect to maps between local coefficient systems over a single base space (this follows for intersection homology with local coefficients just as for ordinary homology with local coefficients).

By the naturality of homology, this chain map induces the desired map of $E^2$
terms of the spectral sequences. Note that, by the formula for the intersection homology of a cone (see \cite{Ki}), this map is the identity for $
*<n-k-1-\bar p(n-k)$ and the
the $0$ map for $*\geq n-k-1-\bar p(n-k)$, since this is true for
each map
$IH_*^{\bar p}(L;\mc{G})\to IH_*^{\bar p}(cL;\mc{G})$ induced by $L\into cL$.
\end{proof}

\begin{proof}[Proof that the right face of the cube commutes] 
We conclude by showing that the right side of the cube \eqref{E: cube} commutes. 
Actually, we begin with the corresponding cube for the $\bar E^1$ terms of the spectral sequence for the intersection homology of the stratified fibration $\hl_s(N,X_k)\to X_k$. The proof is a modified version of that of \cite[XIII.4.7]{Wh}, which is again largely simplified by the fact that the maps are mostly inclusions.

Replacing each term with its direct sum decomposition, we must prove that the following 
diagram commutes:

{\footnotesize\begin{equation*}
\begin{CD}
\oplus_{\beta<\alpha} H_{p-1}(\hat\sigma^{p-1}_{\beta}, \bd \hat\sigma^{p-1}_{\beta };IH^{\bar p}_q(\mc L_{\alpha};\mc{G}| \mc L_{\alpha }) )& 
@>i_{\alpha*}>>
&  \oplus_{\beta} H_{p-1}(\hat\sigma^{p-1}_{\beta}, \bd
\hat\sigma^{p-1}_{\beta};IH^{\bar p}_q(\mc L_{\beta};\mc{G}| \mc L_{\beta}) )\\
@V\lambda' VV && @V\lambda VV\\
\oplus_{\beta<\alpha} IH_{p+q-1}^{\bar p}(\bar \sigma^{p-1}_{\beta}, \bd \bar  \sigma^{p-1}_{\beta};\mc{G}|\bar \sigma^{p-1}_{\beta}) 
&@>j_{\alpha*}>>
&\oplus_{\beta} IH_{p+q-1}^{\bar p}(\bar \sigma^{p-1}_{\beta}, \bd \bar
\sigma^{p-1}_{\beta};\mc{G}|\bar \sigma^{p-1}_{\beta}).\\
\end{CD}
\end{equation*}}

\noindent The top and bottom maps are induced by inclusions and, for $i_{\alpha*}$, the
map
of the local coefficient system determined by the path class in
$\hat\sigma_{\alpha}$
between the basepoints of the generators of
$H_p(\hat\sigma_{\alpha},\bd\hat\sigma_{\alpha})$ and
$H_{p-1}(\hat\sigma_{\beta},\bd\hat\sigma_{\beta})$. The issue is that $\lambda$
is defined using the strong stratified trivialization over each
$\hat \sigma^{p-1}_{\beta}$, while $\lambda'$ is defined using the restriction of the
strong stratified trivialization over $\hat \sigma^p_{\alpha}$. It suffices to prove that this
diagram commutes on each summand. It is clear that if $\beta\neq\beta'$, then the
maps $ H_{p-1}(\hat\sigma^{p-1}_{\beta}, \bd \hat\sigma^{p-1}_{\beta};IH^{\bar
p}_q(\mc L_{\alpha};\mc{G}| \mc L_{\alpha}) )\to IH_{p+q-1}^{\bar p}(\bar
\sigma^{p-1}_{\beta'}, \bd \bar \sigma^{p-1}_{\beta'};\mc{G}|
\sigma^{p-1}_{\beta})$ are each the zero map, so it suffices to show that
\begin{equation*}
\begin{CD}
H_{p-1}(\hat\sigma^{p-1}_{\beta}, \bd \hat\sigma^{p-1}_{\beta};IH^{\bar p}_q(\mc L_{\alpha };\mc{G}| \mc L_{\alpha }) )& 
@>i_{\alpha*}>>
&  H_{p-1}(\hat\sigma^{p-1}_{\beta}, \bd \hat\sigma^{p-1}_{\beta};IH^{\bar p}_q(\mc L_{\beta};\mc{G}| \mc L_{\beta}) )\\
@V\lambda' VV && @V\lambda VV\\
 IH_{p+q-1}^{\bar p}(\bar \sigma^{p-1}_{\beta}, \bd \bar  \sigma^{p-1}_{\beta};\mc{G}| \bar\sigma^{p-1}_{\beta}) 
&@>j_{\alpha*}>>
& IH_{p+q-1}^{\bar p}(\bar \sigma^{p-1}_{\beta}, \bd \bar  \sigma^{p-1}_{\beta};\mc{G}| \bar\sigma^{p-1}_{\beta})\\
\end{CD}
\end{equation*}
commutes. 

The bottom map is a canonical isomorphism, while, again, the top is induced by
the identity of the space $\hat \sigma_{\beta}$ and the
map on coefficients, $u$, determined over the path class in
$\alpha$ between the basepoints of the generators of
$H_p(\hat\sigma_{\alpha},\bd\hat\sigma_{\alpha})$ and
$H_{p-1}(\hat\sigma_{\beta},\bd\hat\sigma_{\beta})$.

So, applying the universal coefficient theorem and the K\"unneth theorem for
intersection homology, we identify $H_{p-1}(\hat\sigma^{p-1}_{\beta}, \bd
\hat\sigma^{p-1}_{\beta};IH^{\bar p}_q(\mc L_{\beta};\mc{G}| \mc L_{\beta}) )$
with
\begin{equation*}
IH^{\bar p}_{p-1}(\hat\sigma^{p-1}_{\beta}\times \mc L_{\beta} , \bd
\hat\sigma^{p-1}_{\beta}\times \mc
L_{\beta};\hat\sigma_{\beta}^{p-1}\times (\mc{G}|\mc{L}_{\beta})),
\end{equation*}
and similarly
for $\alpha$. From the definitions of $\lambda$ and $\lambda'$, we need to show the
commutativity of

{\scriptsize\begin{equation}\label{D: comm}
\begin{CD}
H_{p-1}(\hat\sigma^{p-1}_{\beta}, \bd \hat\sigma^{p-1}_{\beta};IH^{\bar p}_q(\mc L_{\alpha };\mc{G}| \mc L_{\alpha }) )& 
@>i_{\alpha*}>>
&  H_{p-1}(\hat\sigma^{p-1}_{\beta}, \bd \hat\sigma^{p-1}_{\beta};IH^{\bar p}_q(\mc L_{\beta};\mc{G}| \mc L_{\beta}) )\\
@VVV && @VVV\\
IH^{\bar p}_{p-1}(\hat\sigma^{p-1}_{\beta}\times \mc L_{\alpha} , \bd
\hat\sigma^{p-1}_{\beta}\times \mc L_{\alpha};\hat\sigma^{p-1}_{\beta}\times  
(\mc{G}|{\mc{L}_{\alpha}}))& @>(\text{id}\times u)_* >>&
IH^{\bar p}_{p-1}(\hat\sigma^{p-1}_{\beta}\times \mc L_{\beta} , \bd
\hat\sigma^{p-1}_{\beta}\times \mc L_{\beta};\hat\sigma^{p-1}_{\beta}\times
(\mc{G}|{\mc{L}_{\beta}}))\\
@V \mf{h}_{\alpha*}VV && @VV \mf{h}_{\beta*}V\\
IH_{p+q-1}^{\bar p}(\bar \sigma^{p-1}_{\beta}, \bd \bar  \sigma^{p-1}_{\beta};\mc{G}| \bar\sigma^{p-1}_{\beta}) 
&@>j_{\alpha*}>>
& IH_{p+q-1}^{\bar p}(\bar \sigma^{p-1}_{\beta}, \bd \bar  \sigma^{p-1}_{\beta};\mc{G}| \bar\sigma^{p-1}_{\beta}).\\
\end{CD}
\end{equation}}

The bottom vertical maps $\mf{h}_{\alpha*}$ and  $\mf{h}_{\beta*}$ are induced by
$\mf{h}_{\alpha}$ and $\mf{h}_{\beta}$ and covered by the respective coefficient
maps $[\hat \sigma_{\alpha}\times
(\mc{G}|\mc{L}_{\alpha})]|(\hat\sigma_{\beta}\times \mc{L}_{\alpha})\to
(\mf{h}_{\alpha}^*\mc{G})|(\hat\sigma_{\beta}\times \mc{L}_{\alpha})\to \mc{G}$
and $\hat \sigma_{\beta}\times
(\mc{G}|\mc{L}_{\beta})\to
\mf{h}_{\beta}^*\mc{G}\to \mc{G}$. Note that the first map of each of these
compositions is
given by a unique isomorphism compatible with the identity over the base
point. The map $u$ will be the map determined by extending the inclusion of $\mc{L}_{\alpha}$ over a path.

The top half of this diagram commutes by the naturality of the universal
coefficient and K\"unneth theorems. To show that the bottom commutes, we find a
stratum-preserving homotopy between $\mf h_{\alpha}$ and
$\mf h_{\beta}\circ(\text{id}\times
u)$, covered by a coefficient homotopy $\hat\sigma^{p-1}_{\beta}\times
(\mc{G}|{\mc{L}_{\alpha}})\times I\to \mc{G}|\bar\sigma^{p-1}_{\beta}$.  For
this, we will
use a stratified lifting extension problem.

Let $e_{\alpha}$ and $e_{\beta}$ be the base points for the singular 
simplices which generate
$H_p(\hat\sigma_{\alpha},\bd\hat\sigma_{\alpha})$ and
$H_{p-1}(\hat\sigma_{\beta},\bd\hat\sigma_{\beta})$. Given a path $v:I\to \hat
\sigma_{\alpha}$ from $e_{\beta}$ to $e_{\alpha}$ and the projection
$p:\mc{L}_{\alpha}\times I\to I$, let $\td v$ denote a stratum
preserving extension over $vp:\mc{L}_{\alpha}\times I\to B$ of the inclusion map
$\mc L_{\alpha}\into \hl_s(N;X_k)$, and let
$u_t$ be the induced stratum-preserving map from $\mc L_{\alpha}$ to
$\pi^{-1}(v(t))$ (so that $u=u_0$ is the map induced between the fibers over the
endpoints).

We first construct a stratum-preserving homotopy from $\mf 
h_{\alpha}i_{e_{\beta}}$
to $u$, where $ i_{x}:\mc L_{\alpha}\to \hat \sigma_{\alpha}\times \mc L_{\alpha}$ is
the inclusion over $x\in \hat\sigma_{\alpha}$. Consider the diagram
\begin{equation*}
\begin{CD}
\dot I\times I\times \mc L_{\alpha}\cup I\times 1 \times \mc L_{\alpha}  &@>f >>&\bar \sigma_{\alpha}\\
@VVV & & @V\pi VV\\
I\times I\times \mc L_{\alpha}&@>g>>& \hat \sigma_{\alpha},
\end{CD}
\end{equation*}
where $g(t,s,l)=v(s)$ and 
\begin{equation*}
f(t,s,l)=
\begin{cases}
u_s(l) ,&t=0\\
\mf h_{\alpha}i_{v(s)},&t=1\\
\mf h_{\alpha}(l), & s=1.
\end{cases}
\end{equation*}
Since $\mf h_{\alpha}$ is a strong stratified trivialization, $f$ is well-defined,
and it is stratum-preserving by the properties of the component maps. If $K$ is
the product of the identity map on $\mc L_{\alpha}$ with a deformation retraction
of the DR-pair $(I\times I,\dot I\times I\cup I\times 1)$, then $fK_1$ is
stratum-preserving along $I$ and, by Lemma \ref{L: DR ext}, there exists a
stratified solution, $G$, to the lifting extension problem. It is easy to check
that $G$ is stratum-preserving along $I\times I$ (it is stratum-preserving along
one $I$ by the lemma and along the other because it is along $s=1$). In
particular, the restriction of $G|I\times 0\times \mc{L}_{\alpha}$ over
$v(0)=e_{\beta}$, being stratum-preserving
along $I$, gives a stratified homotopy between $\mf h_{\alpha}i_{e_{\beta}}$ and
$u$.
Let us denote this homotopy by $\td u: I\times \mc L_{\alpha}\to \mc L_{\beta}$.

Now consider the diagram
\begin{equation*}
\begin{CD}
\dot I\times\hat\sigma_{\beta}\times \mc L_{\alpha}\cup I\times \{e_{\beta}\}\times \mc L_{\alpha} &@>\eta>>&\bar\sigma_{\beta}\\
@Vi VV&&@V\pi VV\\
I\times \hat\sigma_{\beta}\times \mc L_{\alpha} &@>\tau>>&\hat\sigma_{\beta},
\end{CD}
\end{equation*}
where $\tau(t,s,l)=s$ and
\begin{equation*}
\eta(t,s,l)=
\begin{cases}
\mf h_{\alpha}(s,l), & t=1\\
\mf h_{\beta}(s,u(l)), & t=0\\
\td u(t,l), & s=e_{\beta}.
\end{cases}
\end{equation*}
As $\mf h_{\beta}$ is a strong stratified trivialization and $\td u$ is a
stratified homotopy from $u$ to $h_{\alpha}i_{e_{\beta}}$, $\eta$ is well defined.
Furthermore, by the stratum-preserving properties of the components, $\eta$ is
stratum-preserving with respect to the product filtration of $ I\times
\hat\sigma_{\beta}\times \mc L_{\alpha}$. Therefore, if $K'$ is the product of the
identity map on $\mc L_{\alpha}$ with a deformation retraction of the DR-pair
$(I\times \sigma_{\beta},\dot I\times \sigma_{\beta}\cup e_{\beta}\times I)$ (see
\cite[I.5.2]{Wh}), then $\eta K'$ is stratum-preserving along $I$ and, by Lemma
\ref{L: DR ext}, there exists a stratified solution, $\td \eta$ to the lifting
extension problem. As $\td \eta$ is stratum-preserving along $I$, $\td \eta$ is a
stratum-preserving map, and, in particular, a stratum-preserving homotopy between
$\mf h_{\alpha}$ and $\mf h_{\beta}\circ(\text{id}\times u)$. In addition, these
constructions
can all be lifted in the exact same manner to give covering homotopies of
all of the local coefficient maps  (see also the section on the
construction of
the bundle
of coefficients of a stratified fibration, above). Hence $\mf h_{\alpha}$ and $\mf
h_{\beta}\circ(\text{id}\times u)$, covered by their coefficient maps, induce
the same maps on intersection homology and diagram \eqref{D: comm} commutes.

This completes the proof of the commutativity of the right side of the cube for the
$\bar E^1$ terms. For the $E^1$ terms, 
the commutativity follows from that of the $\bar E^1$ terms by our standard methods
for turning maps and homotopies of stratified fibrations into maps and homotopies
of their mapping cylinders.
In
slightly greater detail: we could build another cube with the $\bar E^1$ terms on
the right, the $E^1$ terms on the left, and with maps between them induced by
the obvious inclusions. The front, back, top, and bottom commute by the naturality
of
inclusions,
the naturality inherent in our definitions of $\lambda$ and $\lambda'$,  the
nature of the construction of the $h_{\gamma}$ from the $\mf{h}_{\gamma}$, and the
same for all of the coefficient covers. We
have just
proven that the left side commutes, hence so does the right.
Also, we could note that the right face of the cube for the $E^1_{*,q}$
terms is
isomorphic to that for the $\bar E^1_{*,q}$ terms if 
$ *<n-k-1-\bar p(n-k)$ and identically zero for $*\geq
n-k-1-\bar(n-k)$ (again, by the formula for the intersection homology of a cone in \cite{Ki}).  

\end{proof}

This, at last,  completes the proof of Theorem \ref{P: Neigh IH}.

\appendix

\section{The stratum-preserving homotopy equivalence of a deformation retract neighborhood of a pure subset with the mapping cylinder of its holink}

In this appendix, we provide a proof of  a theorem found in works of Chapman \cite{Ch79} and Quinn \cite{Q1}. It is necessary to provide some details since the proofs in those works make use of a map that turns out not to be continuous. The proof below was suggested by Bruce Hughes. 

Suppose that $Y$ is a manifold weakly stratified space and that $X$ is a \emph{pure} subset, i.e. $X$ is a closed subspace that is the union of connected  strata $Y_k-Y_{k-1}$. The subsets of interest in  this paper are the connected components of bottom strata of $Y$. Suppose also that $X$ has a neighborhood $N$ in $Y$ and that there exists a nearly stratum-preserving strong deformation retraction $R: N\times I\to N$ that takes $N$ into $X$ (recall that \emph{nearly-stratum preserving} means that $R$ respects the stratification except at time $0$ when all of $N$ is mapped into $X$).  The theorem asserts that under these assumptions, $N$ is stratum-preserving homotopy equivalent to the mapping cylinder of the holink evaluation map $\pi: \hl_s(N,X)\to X$ given by evaluation of paths at time $0$.:

\begin{proposition}
Given a nearly stratum-preserving deformation retraction $R:N\times I\to N$ of $N$ to $X$, there exists a stratum-preserving homotopy equivalence between $N$ and the mapping cylinder $I_{\pi}$  of the holink evaluation $\pi:\hl_s(N,X)\to X$.
\end{proposition}

\begin{proof}
Recall that we assume the \emph{teardrop topology} on the mapping cylinder $I_{\pi}$ of the holink evaluation $\pi: \hl_s(N,X)\to X$, i.e. the topology generated by the basis of open subsets of $\hl_s(N,X)\times (0,1]$ with the product topology and sets of the form $[\pi^{-1}(U)\times (0,\epsilon)]\cup U$, where $U$ is an open set in $X$. 

Also recall (\cite[Lemma 2.4]{Q1} and \cite{Hug}) that given a map $\delta: X\to (0,\infty)$, the inclusion $\hl_s^{\delta}(N,X)\into \hl_s(N,X)$ is a stratum and fiber preserving homotopy equivalence (where $\hl_s^{\delta}(N,X)$ is the subset of paths $w$ in $\hl_s(N,X)$ that are contained entirely within distance $\delta(w(0))$ of their endpoints $w(0)$).  Furthermore, the proof of this fact comes from constructing a deformation retraction of $ \hl_s(N,X)$ into $\hl_s^{\delta}(N,X)$ by shrinking paths along themselves. Similarly, we can note that if $\delta'(x)<\delta(x)$ then $\hl_s^{\delta'}(N,X)\into \hl_s^{\delta}(N,X)$ is also a stratum and fiber preserving homotopy equivalence also given by a deformation retraction that shrinks paths and hence remains entirely in $\hl_s^{\delta}(N,X)$. In fact, just restrict the retraction of the whole space $\hl_s(N,X)$ to $\hl_s^{\delta}(N,X)$. We now define a shrinking map $S: \hl_s(N,X)\times (0,1]\to \hl_s(N,X)$ as follows: on $\hl_s(N,X)\times [1/2,1]$, define $S$ as a deformation retraction that shrinks $\hl_s(N,X)$ into $\hl_s^{1/2}(N,X)$ (scaled to occur on $[1/2,1]$ instead of $[0,1]$); on $[1/3,1/2]$, define $S$ as the composition of $S(\cdot, 1/2)\times \text{id}_{[1/3,1/2]}$ followed by a  deformation retraction from $\hl_s^{1/2}(N,X)$ into $\hl_s^{1/3}(N,X)$ (again scaling $I$ to the interval $[1/3,1/2]$); and so on. So $S$ is stratum and fiber preserving continuous homotopy that shrinks all of the paths in the holink. We can extend $S$ to a map $\bar S: I_{\pi}\to I_{\pi}$ as follows: Let $\bar  S$ be the the identity on the base $X$, and on $\hl_s(N,X)\times (0,1]$ define $\bar S$ by $\bar S(w,s)=(S(w,s),s)$. The map $\bar S$ is continuous: this is clear on $I_{\pi}-X$; for a point $x\in X$, notice that any neighborhood of the form $[\pi^{-1}(U)\times (0,\epsilon)]\cup U$, $U$ an open subset of $X$ with $x\in U$, gets mapped into itself under $\bar S$. 

We can now define the stratum-preserving homotopy inverses $f:N\to I_{\pi}$ and $g:I_{\pi}\to N$. Given a point $y\in N$, let $R_y$ be the retraction path of $y$ under the nearly stratum-preserving deformation retraction $R$. Define $f$ by $f(y)=y$ if $y\in X$ and $f(y)=(R_y, d(y, X))$ if $y\in N-X$. Here $d(y, X)$ is the distance from $y$ to $X$. Define $g$ by $g(y)=y$ for $y\in X$ and for $(w,s)\in \hl_s(N,X)\times (0,1]$ by $g(w,s)=S(w,s)(s)$; this is the evaluation of the path $S(w,s)$ at the time $s$. We first show that $f$ and $g$ are continuous, then that they are in fact stratum-preserving homotopy inverses.

For $f$, continuity is clear except at points in $X$. So choose $x\in X$, and suppose that we are given an open set $V\subset I_{\pi}$ of the form $[\pi^{-1}(U)\times (0,\epsilon)]\cup U$ such that $U$ is an open set in $X$ that contains $f(x)=x$. We need to find a neighborhood of $x$ in $N$ that maps into $V$. Let $\delta=\text{min}(\epsilon, d(x, X-U))$, where the distance $d$ is that of  $N$ restricted to $X$. Let $B(x, \delta)$ be the open neighborhood of $x$ in $N$ of radius $\delta$. Consider the set $R^{-1}(B)\subset N\times I$. Note that this open set includes $B\cap X$, and so by standard general topology (see, e.g., \cite{MK2}), it includes a set of the form $W\times I$ where $W$ is an open neighborhood of $B\cap X$ in $N$. Hence, in particular, $W$ is an open neighborhood of $x$. So let us consider $f(w)$ for $w\in W$. If $w\in W\cap X$, then $w\in B\cap X\subset U$ so that $f(w)\in V$. If $w\in W-X\cap W$, then $f(w)=(R_w, d(w,X))$. The path $R_w$ is the image of $R(w, \cdot)$, but by construction, this path lies entirely within $B$. Hence the $0$ endpoint must lie in $B\cap X\subset U$, and $R_w\subset \pi^{-1}(U)$. Also, since $w =R(w,1)\in B$, $d(w,X)<\epsilon$.  This shows that $f(w)\subset V$. We conclude that $f(W)\subset V$, and since our choice of $V$ was arbitrary, $f$ is continuous.

Next, we show that $g$ is continuous. Again, continuity is clear at points not in $X\subset I_{\pi}$. So let $x\in X$, and consider a neighborhood of $x\in N$. Since $N$ has the metric topology, we can assume that this neighborhood has the form $B=B(x,\epsilon)$, the ball of radius $\epsilon$ about $x$ in $N$. So let $U$ be the open ball about $x$ in $X$ of radius $\epsilon/2$, and choose $n$ so that $1/n<\epsilon/2$. Consider the subset of $I_{\pi}$ given by $V=[\pi^{-1}(U)\cup (0, 1/n)]\cup U$. We will see that $g(V)\subset B$. If $x\subset U$, then $g(x)=x\in U\subset B\cap X$. For points of the form $(w,s)\in \hl_s(N,X)\times (0,1]$, $g(w,s)=S(w,s)(s)$. Since $s<1/n<\epsilon/2$ by assumption, it follows from the construction of $S$ that $S(w,s)$ is a path whose image lies within $\epsilon/2$ of its endpoint in $U$. So $d(g(w,s),x)\leq d(g(w,s), w(0))+d(w(0),x)<\epsilon/2+\epsilon/2=\epsilon$, the latter $\epsilon/2$ occurring because we have assume $w(0)\subset U$ (since $w\in \pi^{-1}(U)$) and $U$ has radius $\epsilon/2$ from $x$. Thus $g$ is continuous. 

Now let us show that $gf$ is stratum-preserving homotopy-equivalent to the identity. Again, if $y\in X$, $gf(y)=y$, but if $y\in N-X$, $gf(y)=g(R_y, d(y, X))=S(R_y,d(y,X))(d(y,X))$. Since $S$ retracts paths along themselves, $gf(y)=R_y(s')=R(y,s')$, where $s'$ is a time that can be chosen to depend continuously on $y$. In fact, $S(w,s)(t)$ is given by $w(\xi(w,s,t))$, where $\xi$ is a continuous map  that we can think of as a reparametrization from the $t$ factor $[0,1]$ into (but not necessarily onto) $[0,1]$ which takes $0$ to $0$, $(0,1]$ into $(0,1]$, and depends continuously on the parameters $w$ and $s$. Then given $y$,  $s'=R_y(\xi(R_y, d(y,X), d(y,X))$. The homotopy from $gf$ to the identity is then given by $(y,t)\to R(y,t+(1-t)s')$: at $t=0$, $y\to gf(y)$ and at $t=1$, $y\to R(y,1)=y$. This homotopy is stratum-preserving because $R$ is. Again, there are no problems seeing the continuity at points in $N\times I-X\times I$. For points $(x,t)\in X\times I$, we can emulate the arguments whereby we showed that $f$ was continuous by looking only at points in a small neighborhood of $x$ that also have retraction paths contained entirely in a small neighborhood of $x$. 

Finally, we must show that $fg$ is stratum-preserving homotopy-equivalent to the identity. Again, let us first calculate $gf$. For $x\in X\subset I_{\pi}$, $fg(x)=x$, while for $(w,s)\in \hl_s(N,X)\times (0,1]\subset I_{\pi}$, $fg(w,s)=f(S(w,s)(s))=(R_{S(w,s)(s)}, d(S(w,s)(s), X))$. Again, let us choose $s'$ so that $S(w,s)(s)=w(s')$ and note that $s'$ can be chosen to depend continuously on $w$ and $s$. 
We now define a homotopy $h: I_{\pi}\times I\to I_{\pi}$: let $h|X\times I$ be the projection onto $X$,  and let the first coordinate   of the image of $\hl_s(N,X)\times (0,1]$ be determined by the paths given as follows (where $p_i$ is the projection of $\hl_s(N,X)\times (0,1]$ onto its $i$th factor):

\begin{equation}\label{E: hl he}
p_1h((w,s), u)(t)=
	\begin{cases}
	R(w(t),1),
		\hfill u=1\\
	R(w(ut+(1-u)s'), \frac{t}{t-u+1}), 
            	\hfill 0\leq t\leq s'u, 0\leq u< 1 \\
	R(w(ut+(1-u)s'), \frac{s'u}{s'u-u+1}(1-\frac{t-s'u}{s'-s'u})+(u+(1-u)s')(\frac{t-s'u}{s'-s'u}),\\
		\hfill s'u\leq t\leq s', 0\leq u<1\\
	R(w(tu+(1-u)s'), u+(1-u)t), 
		\hfill s'\leq t\leq 1, 0\leq u< 1\\
	\end{cases}
\end{equation}
For the second coordinate, we have $p_2h((w,s),u)=us+(1-u)d(S(w,s)(s), X)$. 

The reason these formulas need be so complicated is that we must homotop the paths in such a manner that all paths are in $X$ at time $t=0$ (i.e. to keep all paths in the holink). For each path $w$, the first coordinate of this homotopy takes it to the path of the retraction of the point $w(s')$, sweeping through the retractions of the points $w(t)$, $0\leq t\leq s'$ along the way. 

Note that when $u=1$, we have $p_1h((w,s),1)(t)=R(w(t),1)=w(t)$, and for $u=0$, we have $p_1h((w,s),0)(t)=R (w(s'), t)=R(S(w,s)(s),t)=R_{S(w,s)(s)}(t)$. Also, $p_1h((w,s),u)(0)=R( w((1-u)s'),0)\in X$, and so long as $s>0$, the second input into $R$ in the above formulae will be $>0$ so long as $t>0$.  So, for all $u$ the first coordinate gives us a legitimate element of $\hl(N,X)$. That our homotopy is stratum-preserving follows form the fact that $R$ is stratum-preserving and that paths in $\hl_s(N,X)$ are nearly stratum-preserving. 

To show the continuity of $h$, we divide into two cases, for points in $X\times I$ and for points in $\hl_s(N,X)\times (0,1]\times I$. The latter points map into $\hl_s(N,X)\times (0,1]$ and so to show continuity here we need only check continuity of the component functions (or equivalently, of $p_1h$ and $p_2h$). The continuity of $p_2h$ on $\hl_s(N,X)\times (0,1]\times I$ is clear from the definition.
Let us show the continuity of $p_1h$ on $\hl_s(N,X)\times (0,1]\times I$. By general topology, it suffices to show that $E(p_1h\times \text{id}_I): \hl_s(N,X)\times (0,1]\times I\times I\to N$ is continuous, where $E$ is the holink evaluation $(w,t)\to w(t)$. Observe that for the chosen  domain $s>0$ and so $s'>0$. Looking at equation \eqref{E: hl he}, we see that this map is clearly continuous except possibly for points where $u=1$ and $t$ is equal to $s'$ or $0$. So let us directly check continuity at such points. 

The potential difficulty at points $(w,s,s',1)$ is, in fact, illusory. The above equations have been written to make clear where linear interpolations are occurring, and a rewriting of the third line of the equation by simplifying the algebra will illustrate that the denominator occurs in a neighborhood of $(s')^2$ when $u$ is near $1$. It is then easy to check that $p_1h$ is continuous in the neighborhoods of such points.  

Consider then a point $(w_0,s_0,1,0)$. Near this point, the governing equations will be $R(w(t),1)$ for $u=1$ and $R(w(ut+(1-u)s'), \frac{t}{t-u+1})$ for $ u< 1$.
By definition, $h((w_0,s_0),1)(0)=w_0(0)\in X\subset N$. Let $B$ be an $\epsilon$ neighborhood of $w_0(0)$ in $N$. We find an open neighborhood of $(w_0,s_0,1,0)$ that maps into $B$. Without constructing such a neighborhood explicitly we indicate that one is easy to obtain using the product topology on $\hl_s(N,X)\times (0,1]\times I\times I$. Let us briefly keep fixed $w_0$ and  $s_0$ (and hence  $s'_0=s'(w_0,s_0)$). Now, as observed in the proof of continuity of $f$, given a $\delta$ (say $\epsilon/2$), $x$ has a neighborhood $U$ in $N$ such that any retraction path under $R$ that starts in this neighborhood remains $\delta$ close to $x$. If we take $t$ sufficiently small, then $w_0(t)$ will be in this neighborhood and so will $ R(w_0(ut+(1-u)s'_0), \frac{t}{t-u+1})$ because these points will lie in the retractions paths of the points on $w_0$ in the neighborhood. Thus we obtain continuity on the restriction to fixed $w_0$ and $s_0$. But now observe that since we use the compact-open topology on $\hl_s(N,X)$, we can choose a neighborhood of $w_0$ in $\hl_s(N,X)$ so that all paths in this neighborhood are very close to the path $w_0$. We can also choose a small interval about $s_0$ and since $s'$ depends continuously on $w$ and $s$, we can choose these neighborhoods so that $s'$ also varies only slightly. Hence it is clear that for sufficiently refined choices of neighborhoods of $w_0$ and $s_0$ and for the same choices of sufficiently small neighborhoods of $t=0$ and $u=1$, the image points under  
$ R(w(ut+(1-u)s'), \frac{t}{t-u+1})$ or $R(w(t),1)$ will remain close to $w(0)$. In particular, with $u$ sufficiently close to $1$ and $t$ sufficiently close to $0$, all the $ w(ut+(1-u)s')$ will remain within our  neighborhood $U$ of $w_0(0)$ and so the retraction paths utilized will all remain within $\epsilon/2$ of $w_0(0)$.

To prove that $h$ is a continuous homotopy on all of $I_{\pi}\times I$, it remains only to check continuity at points of the form $(x,t)\in X\times I\subset I_{\pi}\times I$.
Let $V$ be a basis element of $I_{\pi}$ that contains $h(x,t)=x$. Then $V$ is of the form $[\pi^{-1}(U)\times (0,\epsilon)]\cup U$, and we need to find an open neighborhood $W$ of $(x,t)$ in $I_{\pi}\times I$ that maps into $V$ under $h$. Recall (Lemma \ref{L: prod/cyl}) that the topology of $I_{\pi}\times I$ is the same as that given if we consider this set as the mapping cylinder of $\pi\times \text{id}_I$. So we will try to find a $W$ of the form $[(\pi\times \text{id}_I)^{-1}(Z)\times (0,\delta)]\cup Z$, where $Z$ is an open neighborhood of $(x,t)$ in $X\times I$. 

Let $D$ be the open ball about $x$ in $N$ of radius $\rho'<\text{min}(\epsilon/2,d(x, X-U))$, where the distance $d$ is that obtained from $N$ by restricting to $X$. Consider $R^{-1}(D)\subset N\times I$, and let $C\times I$ be an open neighborhood of $(D\cap X)\times I$ within $R^{-1}(D)$ (again, such a neighborhood exists by standard topology). Now let $B$ be a  neighborhood of $x $ in $X$ of radius $\rho<\rho'/2$, and choose $\delta'$ so that $0<\delta'< d(B, N-C)$ (note, of course, that it is possible to choose such a $\delta'$ by construction). 

Let us choose $n$ such that $1/n<\delta'$ and $1/n<\epsilon/2$, and choose $\delta<\text{min}( \rho, 1/n)$. Let $B$ be as above, the open ball about $x$ in $X$ of radius $\rho$, and let $Z=B\times I$.  Then $Z\subset U\times I$, and with this choice of $Z$ and $\delta$, let $W$ be defined as above. We will show that $h(W)\subset V$. If $(x,u)\in Z\times I\subset X\times I$, then $h(x,t)=x\subset B\subset U\subset V$. So let us consider a point in $W$ of the form $((w,s),u)$. We need to show that $h((w,s),u)\in \pi^{-1}(U)\times (0,\epsilon)$.  Since $p_2h((w,s),u)=us+(1-u)d(S(w,s)(s), X)$, $p_2h$ lies in $(0,1]$ between $s$ and $d(S(w,s)(s), X)$. By assumption, $s<\delta<\epsilon/2$, and $d(S(w,s)(s), X)<\epsilon/2$ since $S(w,s)$ is a path contained within $\epsilon/2$ of $w(0)\in X$ by the construction of the map $S$. Thus $p_2h((w,s),u)<\epsilon$. So finally, we need only check that $p_1h((w,s),u)(0)\subset U$ to conclude that $h(W)\subset V$. As we observed above, $p_1h((w,s),u)(0)=R( w((1-u)s'),0)$, and as $u$ varies these are simply the end points of the retraction paths of the points on $w$ between $0$ and $s'$.  But then the images of $w$ in this range must also be images of $S(w,s)$ and each such point has a distance $<\delta'$ from $x\subset B$ and so lies in $C$. But then, again by construction, each of these points has a retraction path which lies entirely within $D$ and, in particular then, has its endpoint in $D\cap X \subset U$. At last, we conclude that $h(W)\subset V$ and so $h$ is continuous.

\end{proof}

\bibliographystyle{amsplain}
\bibliography{bib}

Several diagrams in this paper were typeset using the \TeX\, commutative
diagrams package by Paul Taylor.

\end{document}